\documentclass[a4paper,11pt,twoside,english]{article}
\usepackage{babel}
\usepackage{latexsym,amsfonts}
\usepackage{amsthm}
\usepackage{amsmath}
\usepackage{mathtools}
\usepackage{paralist}
\usepackage{ifthen}
\usepackage{moredefs}
\usepackage{todonotes}
\usepackage{tikz}
\usepackage{booktabs}
\usepackage{float}

\usepackage{caption}
\providecommand{\CAP}[1][-30pt]{
  \captionsetup{font=small,name=Fig.,justification=raggedleft,singlelinecheck=false,skip=#1}
  \caption{}
  }

\usepackage{lastpage}
\usepackage{fancyhdr}
\fancypagestyle{plain}{

  \fancyhead[L,C]{}
  \fancyhead[R]{\today}
  \fancyfoot[R,L]{}
  \fancyfoot[C]{\thepage \ of \pageref{LastPage}}
  }

\title{K\"othe-Bochner spaces and some geometric properties related to rotundity and smoothness}
\author{Jan-David Hardtke}
\date{}

\providecommand{\shortversion}[1]{#1}
\providecommand{\longversion}[2]{#2 {\rm(#1 in short)}}

\providecommand{\acronym}[2]{
  \newboolean{acronym#1}
  \setboolean{acronym#1}{true}
  \expandafter\providecommand\expandafter{\csname acronym#1\endcsname}{#2}
  }

\providecommand{\ac}[1]{\ifthenelse{\boolean{acronym#1}}
{\longversion{#1}{\csname acronym#1\endcsname}\Global\ToggleBoolean{acronym#1}}
{\shortversion{#1}}}

\providecommand{\ifif}{iff }
\providecommand{\sm}{\setminus}
\providecommand{\ssq}{\subseteq}

\providecommand{\N}{\ensuremath{\mathbb{N}}}

\providecommand{\R}{\ensuremath{\mathbb{R}}}

\providecommand{\A}{\ensuremath{\mathcal{A}}}

\providecommand{\eps}{\ensuremath{\varepsilon}}

\providecommand{\keywords}[1]{
  {\let\thefootnote=\relax
  \footnote{{\em Keywords}: #1}}
  \addtocounter{footnote}{-1}
  }

\providecommand{\AMS}[1]{
  {\let\thefootnote=\relax
  \footnote{{\em AMS Subject Classification} (2010): #1}}
  \addtocounter{footnote}{-1}
  }

\providecommand{\address}{
  {\sc \noindent Department of Mathematics \\
  Freie Universit\"at Berlin \\
  Arnimallee 6, 14195 berlin \\
  Germany \\}
  }

\DeclarePairedDelimiter{\set}{\lbrace}{\rbrace}
\DeclarePairedDelimiter{\paren}{\lparen}{\rparen}

\DeclarePairedDelimiter{\abs}{\lvert}{\rvert}
\DeclarePairedDelimiter{\norm}{\lVert}{\rVert}

\theoremstyle{definition}
\newtheorem{definition}{Definition}[section]
\newtheorem*{definition*}{Definition}

\newtheorem*{example*}{Example}

\newtheorem*{remark*}{Remark}

\newtheorem*{problem*}{Problem}

\newtheorem*{question*}{Question}

\newtheorem*{conjecture*}{Conjecture}
\theoremstyle{remark}

\newtheorem*{claim*}{Claim}

\newtheorem*{fact*}{Fact}
\theoremstyle{plain}
\newtheorem{lemma}[definition]{Lemma}
\newtheorem*{lemma*}{Lemma}
\newtheorem{proposition}[definition]{Proposition}
\newtheorem*{proposition*}{Proposition}
\newtheorem{theorem}[definition]{Theorem}
\newtheorem*{theorem*}{Theorem}
\newtheorem{corollary}[definition]{Corollary}
\newtheorem*{corolary*}{Corollary}

\newenvironment{Proof}[1][\proofname]{\begin{proof}[#1] \setlength{\parindent}{0pt}}{\end{proof}}
\newenvironment{Abstract}{\centering\begin{minipage}{0.8\textwidth} \noindent \small {\sc Abstract.}}{\end{minipage}\par}

\usepackage{color}
\definecolor{darkgreen}{rgb}{0,0.5,0}

\numberwithin{equation}{section}
\newtagform{colored}[\color{blue}]{\color{blue}(}{\color{blue})}
\usetagform{colored}

\usepackage[colorlinks,linkcolor=blue,citecolor=red,urlcolor=darkgreen]{hyperref}
\providecommand{\email}{{\it E-mail address:} \href{mailto:hardtke@math.fu-berlin.de}{\tt hardtke@math.fu-berlin.de}}
\providecommand{\mr}[1]{\href{http://www.ams.org/mathscinet-getitem?mr=#1}{MR#1}}

\usepackage{amsrefs}

\begin{document}

\maketitle

\begin{Abstract}
Kadets et al. (cf.\,\cite{kadets}) introduced the notions of acs, luacs and uacs spaces, which form 
common generalisations of well-known rotundity and smoothness properties of Banach spaces. In the 
preprint \cite{hardtke} the author introduced some further related notions and investigated the behaviour 
of these geometric properties under the formation of absolute sums. This paper is in a sense a continuation
of \cite{hardtke}. Here we will study the behaviour of said properties under the formation of K\"othe-Bochner 
spaces, thereby generalising some results of G. Sirotkin from \cite{sirotkin} on the acs, luacs and uacs
properties of $L^p$-Bochner spaces.
\end{Abstract}
\AMS{46B20}

\section{Introduction}\label{sec:intro}
\acronym{UR}{uniformly rotund}
\acronym{LUR}{locally uniformly rotund}
\acronym{WLUR}{weakly locally uniformly rotund}
\acronym{WUR}{weakly uniformly rotund}
\acronym{R}{rotund}
\acronym{S}{smooth}
\acronym{FS}{Fr\'echet-smooth}
\acronym{UG}{uniformly G\^ateaux-differentiable}
\acronym{US}{uniformly smooth}
\acronym{MLUR}{midpoint locally uniformly rotund}
\acronym{WMLUR}{weakly midpoint locally uniformly rotund}
\acronym{acs}{alternatively convex or smooth}
\acronym{luacs}{locally uniformly alternatively convex or smooth}
\acronym{uacs}{uniformly alternatively convex or smooth}
\acronym{sluacs}{strongly locally uniformly alternatively convex or smooth}
\acronym{wuacs}{weakly uniformly alternatively convex or smooth}
\acronym{msluacs}{midpoint strongly locally uniformly alternatively convex or smooth}
\acronym{mluacs}{midpoint locally uniformly alternatively convex or smooth}

We begin with some notation and definitions. Throughout this paper, $X$ denotes a real Banach space, $X^*$ its dual, 
$B_X$ its unit ball and $S_X$ its unit sphere.\par
In the next definition, we summarise the most important rotundity properties.
\begin{definition}\label{def:rotundity}
A Banach space $X$ is called
  \begin{enumerate}[(i)]
    \item {\em\ac{R}} if for any two elements $x,y\in S_X$ the equality $\norm*{x+y}=2$ implies $x=y$,
    \item {\em\ac{LUR}} if for every $x\in S_X$ the implication
    \begin{equation*}
    \norm*{x_n+x}\to 2 \ \Rightarrow \ \norm*{x_n-x}\to 0
    \end{equation*}
    holds for every sequence $(x_n)_{n\in \N}$ in $S_X$,
    \item {\em\ac{WLUR}} if for every $x\in S_X$ and every sequence $(x_n)_{n\in \N}$ in $S_X$ we have
    \begin{equation*}
    \norm*{x_n+x}\to 2 \ \Rightarrow \ x_n\to x \ \mathrm{weakly},
    \end{equation*}
    \item {\em\ac{UR}} if for any two sequences $(x_n)_{n\in \N}$ and $(y_n)_{n\in \N}$ in $S_X$ the implication 
    \begin{equation*}
    \norm*{x_n+y_n}\to 2 \ \Rightarrow \ \norm*{x_n-y_n}\to 0
    \end{equation*}
    holds,
    \item {\em\ac{WUR}} if for any two sequences $(x_n)_{n\in \N}$ and $(y_n)_{n\in \N}$ the following implication holds
    \begin{equation*}
    \norm*{x_n+y_n}\to 2 \ \Rightarrow \ x_n-y_n\to 0 \ \mathrm{weakly}.
    \end{equation*}
  \end{enumerate}
\end{definition}
\noindent The chart below shows the obvious implications between these notions. No other implications are valid in general
(see the examples in \cite{smith}). Note, however, that all these notions coincide in finite-dimensional spaces, by the
compactness of $B_X$.

\begin{figure}[H]
\begin{center}
  \begin{tikzpicture}
  \node (UR) at (-2,0) {UR};
  \node (WUR) at (0,1) {WUR};
  \node (LUR) at (0,-1) {LUR};
  \node (WLUR) at (2,0) {WLUR};
  \node (R) at (4,0) {R};
  \draw[->] (UR)--(LUR);
  \draw[->] (UR)--(WUR);
  \draw[->] (LUR)--(WLUR);
  \draw[->] (WUR)--(WLUR);
  \draw[->] (WLUR)--(R);
  \end{tikzpicture}
\end{center}
\CAP\label{fig:1}
\end{figure}

The modulus of convexity of the space $X$ is defined by 
\begin{equation*}
\delta_X(\eps)=\inf\set*{1-1/2\norm*{x+y}:x,y\in B_X \ \mathrm{and} \ \norm*{x-y}\geq\eps}
\end{equation*}
for every $\eps$ in the interval $]0,2]$. Then $X$ is \ac{UR} \ifif $\delta_X(\eps)>0$ for all $0<\eps\leq2$.\par
For the local version one defines
\begin{equation*}
\delta_X(x,\eps)=\inf\set*{1-1/2\norm*{x+y}:y\in B_X \ \mathrm{and} \ \norm*{x-y}\geq\eps}
\end{equation*}
for every $x\in S_X$ and each $\eps\in ]0,2]$. Then $X$ is \ac{LUR} \ifif $\delta_X(x,\eps)>0$ for all $x\in S_X$ and 
all $0<\eps\leq2$.\par
Let us also recall some notions of smoothness. The space $X$ is called {\em\ac{S}} if its norm is G\^ateaux-differentiable 
at every non-zero point (equivalently at every point of $S_X$), which is the case \ifif for every $x\in S_X$ there is a 
unique functional $x^*\in S_{X^*}$ with $x^*(x)=1$ (cf. \cite{fabian}*{Lemma 8.4 (ii)}). $X$ is called {\em\ac{FS}} if the norm is 
Fr\'echt-differentiable at every non-zero point. The norm of the space $X$ is said to be {\em\ac{UG}} if for each $y\in S_X$
the limit $\lim_{\tau\to 0}\paren*{\norm{x+\tau y}-1}/\tau$ exists uniformly in $x\in S_X$. Finally, $X$ is called {\em\ac{US}} 
if $\lim_{\tau \to 0}\rho_X(\tau)/\tau=0$, where $\rho_X$ denotes the modulus of smoothness of $X$ defined by 
$\rho_X(\tau)=\sup\set*{1/2(\norm*{x+\tau y}+\norm*{x-\tau y}-2):x,y\in S_X}$ for every $\tau>0$.\par
In \cite{kadets} the following notions were introduced (in connection with the so called Anti-Daugavet property).
\begin{definition}\label{def:acs-luacs-uacs}
A Banach space $X$ is called
  \begin{enumerate}[(i)]
  \item {\em\ac{acs}} if for every $x,y\in S_X$ with $\norm*{x+y}=2$ and every $x^*\in S_{X^*}$ with $x^*(x)=1$ we have $x^*(y)=1$ as well,
  \item {\em\ac{luacs}} if for every $x\in S_X$, every sequence $(x_n)_{n\in \N}$ in $S_X$ and every functional $x^*\in S_{X^*}$ we have
  \begin{equation*}
  \norm*{x_n+x}\to 2 \ \mathrm{and} \ x^*(x_n)\to 1 \ \Rightarrow \ x^*(x)=1,
  \end{equation*}
  \item {\em\ac{uacs}} if for all sequences $(x_n)_{n\in \N}$, $(y_n)_{n\in \N}$ in $S_X$ and $(x_n^*)_{n\in \N}$ in $S_{X^*}$ we have
  \begin{equation*}
  \norm*{x_n+y_n}\to 2 \ \mathrm{and} \ x_n^*(x_n)\to 1 \ \Rightarrow \ x_n^*(y_n)\to 1.
  \end{equation*}
  \end{enumerate}
\end{definition}
The author introduced the following related notions in \cite{hardtke}.
\begin{definition}\label{def:sluacs-wuacs}
A Banach space $X$ is called 
  \begin{enumerate}[(i)]
  \item {\em\ac{sluacs}} if for every $x\in S_X$ and all sequences $(x_n)_{n\in \N}$ in $S_X$ and $(x_n^*)_{n\in \N}$ in $S_{X^*}$ we have
  \begin{equation*}
  \norm*{x_n+x}\to 2 \ \mathrm{and} \ x_n^*(x_n)\to 1 \ \Rightarrow \ x_n^*(x)\to 1,
  \end{equation*}
  \item {\em\ac{wuacs}} if for any two sequences $(x_n)_{n\in \N}$, $(y_n)_{n\in \N}$ in $S_X$ and every functional $x^*\in S_{X^*}$ we have
  \begin{equation*}
  \norm*{x_n+y_n}\to 2 \ \mathrm{and} \ x^*(x_n)\to 1 \ \Rightarrow \ x^*(y_n)\to 1.
  \end{equation*}
  \end{enumerate}
\end{definition}
The obvious implication between the \ac{acs} properties and the rotundity properties are indicated in the following chart.
No other implications are generally valid (see the examples in \cite{hardtke}), but note again that the properties \ac{acs},
\ac{luacs}, \ac{sluacs}, \ac{wuacs} and \ac{uacs} coincide in finite-dimensional spaces, by compactness.

\begin{figure}[H]
\begin{center}
\begin{tikzpicture}
  \node (UR) at (-2,0) {UR};
  \node (WUR) at (0,1) {WUR};
  \node (LUR) at (0,-1) {LUR};
  \node (WLUR) at (2,0) {WLUR};
  \node (R) at (4,0) {R};
  \node (uacs) at (-2,-1) {uacs};
  \node (wuacs) at (0,0) {wuacs};
  \node (sluacs) at (0,-2) {sluacs};
  \node (luacs) at (2,-1) {luacs};
  \node (acs) at (4,-1) {acs};
  \draw[->] (UR)--(LUR);
  \draw[->] (UR)--(WUR);
  \draw[->] (LUR)--(WLUR);
  \draw[->] (WUR)--(WLUR);
  \draw[->] (WLUR)--(R);
  \draw[->] (uacs)--(sluacs);
  \draw[->] (uacs)--(wuacs);
  \draw[->] (sluacs)--(luacs);
  \draw[->] (wuacs)--(luacs);
  \draw[->] (luacs)--(acs);
  \draw[->] (UR)--(uacs);
  \draw[->] (LUR)--(sluacs);
  \draw[->] (WUR)--(wuacs);
  \draw[->] (WLUR)--(luacs);
  \draw[->] (R)--(acs);
\end{tikzpicture}
\end{center}
\CAP\label{fig:2}
\end{figure}

The connection between some of the \ac{acs} properties to smoothness properties is illustrated in the diagram below.

\begin{figure}[H]
\begin{center}
\begin{tikzpicture}
  \node (US) at (-2,1) {US};
  \node (UG) at (0,1) {UG};
  \node (S) at (2,1) {S};
  \node (uacs) at (-2,-0.5) {uacs};
  \node (sluacs) at (0,-0.5) {sluacs};
  \node (acs) at (2,-0.5) {acs};
  \draw[->] (US)--(UG);
  \draw[->] (UG)--(S);
  \draw[->] (uacs)--(sluacs);
  \draw[->] (sluacs)--(acs);
  \draw[->] (US)--(uacs);
  \draw[->] (UG)--(sluacs);
  \draw[->] (S)--(acs);
\end{tikzpicture}
\end{center}
\CAP\label{fig:3}
\end{figure}

Let us mention that if we replace the condition $x_n^*(x_n)\to 1$ by $x_n^*(x_n)=1$ for every $n\in \N$ in the
definitions of the properties \ac{uacs} resp. \ac{sluacs} we still obtain the same classes of spaces. For \ac{uacs} spaces
this was first proved by G. Sirotkin in \cite{sirotkin} using the fact that \ac{uacs} spaces are reflexive (see below). For
\ac{sluacs} spaces this characterisation can be proved by means of the Bishop-Phelps-Bollob\'as-theorem (see \cite{hardtke}*{Proposition 2.1}).\par
This characterisation enables us to define the following `\ac{uacs}-modulus' of a given Banach space
(cf.\,\cite{hardtke}*{Definition 1.4}).
\begin{definition}\label{def:uacs-mod}
For a Banach space $X$ we define
  \begin{align*}
  & D_X(\eps)=\set*{(x,y)\in S_X\times S_X:\exists x^*\in S_{X^*} \ x^*(x)=1 \ \mathrm{and} \ x^*(y)\leq 1-\eps} \\
  & \mathrm{and} \ \delta_{\mathrm{uacs}}^X(\eps)=\inf\set*{1-\norm*{\frac{x+y}{2}}:(x,y)\in D_X(\eps)} \ \forall \eps\in ]0,2].
  \end{align*}
\end{definition}
\noindent Then $X$ is \ac{uacs} \ifif $\delta_{\mathrm{uacs}}^X(\eps)>0$ for every $\eps\in ]0,2]$ and we clearly have 
$\delta_X(\eps)\leq\delta_{\mathrm{uacs}}^X(\eps)$ for each $\eps\in ]0,2]$.\par
The above characterisation shows that the class of \ac{uacs} spaces coincides with the class of $U$-spaces introduced by 
Lau in \cite{lau} and our modulus $\delta_{\mathrm{uacs}}^X$ is the same as the modulus of $u$-convexity from \cite{gao1}. 
Also, the notion of $u$-spaces which was introduced in \cite{dhompongsa} coincides with the notion of \ac{acs} spaces.\par
Recall that a Banach space $X$ is said to be uniformly non-square if there is some $\delta>0$ such that for all $x,y\in B_X$
we have $\norm*{x+y}\leq2(1-\delta)$ or $\norm*{x-y}\leq2(1-\delta)$. It is easily seen that \ac{uacs} spaces are uniformly 
non-square and hence by a well-known theorem of James (cf. \cite{beauzamy}*{p.261}) they are superreflexive, as was observed
in \cite{kadets}*{Lemma 4.4}. For a proof of the superreflexivity of \ac{uacs} spaces that does not rely on James' result on
uniformly non-square spaces, see \cite{hardtke}*{Proposition 2.8}.\par
Let us also restate here the following auxiliary result \cite{hardtke}*{Lemma 2.30} (it is the generalisation of \cite{abramovich}*{Lemma 2.1}
to sequences, with a completely analogous proof).
\begin{lemma}\label{aux lemma}
Let $(x_n)_{n\in \N}$ and $(y_n)_{n\in \N}$ be sequences in the (real or complex) normed 
space $X$ such that $\norm*{x_n+y_n}-\norm*{x_n}-\norm*{y_n}\to 0$.\par
Then for any two bounded sequences $(\alpha_n)_{n\in \N}$, $(\beta_n)_{n\in \N}$ 
of non-negative real numbers we also have 
$\norm*{\alpha_nx_n+\beta_ny_n}-\alpha_n\norm*{x_n}-\beta_n\norm*{y_n}\to 0$.
\end{lemma}
Finally, we will need two more definitions from \cite{hardtke}.
\begin{definition}\label{def:luacs+sluacs+}
A Banach space $X$ is called 
  \begin{enumerate}[(i)]
  \item a $\mathrm{luacs}^+$ space if for every $x\in S_X$, every sequence $(x_n)_{n\in \N}$ in $S_X$
  with $\norm*{x_n+x}\to 2$ and all $x^*\in S_{X^*}$ we have
  \begin{equation*}
  x^*(x_n)\to 1 \ \iff \ x^*(x)=1,
  \end{equation*}
  \item a $\mathrm{sluacs}^+$ space if for every $x\in S_X$, every sequence $(x_n)_{n\in \N}$ in $S_X$
  with $\norm*{x_n+x}\to 2$ and all sequences $(x_n^*)_{n\in \N}$ in $S_{X^*}$ we have
  \begin{equation*}
  x_n^*(x_n)\to 1 \ \iff \ x_n^*(x)\to 1.
  \end{equation*}
  \end{enumerate}
\end{definition}
Obviously, every \ac{WLUR} space is $\mathrm{luacs}^+$ and every \ac{LUR} space is $\mathrm{sluacs}^+$.\par
In the next section we will recall some facts on K\"othe-Bochner spaces.\par

\section{Preliminaries on K\"othe-Bochner spaces}\label{sec:prelim koethe bochner}
If not otherwise stated, $(S,\A,\mu)$ will denote a complete, $\sigma$-finite measure space. For $A\in \A$ we denote by $\chi_A$
the characteristic function of $A$.\par
A K\"othe function space over $(S,\A,\mu)$ is a Banach space $(E,\norm{\cdot}_E)$ of real-valued measurable\footnote{i.\,e. $\A$-Borel-measurable} 
functions on $S$ modulo equality $\mu$-almost everywhere\footnote{We will henceforth abbreviate this by $\mu$-a.\,e. or simply a.\,e. 
if $\mu$ is tacitly understood.} such that 
\begin{enumerate}[(i)]
\item $\chi_A\in E$ for every $A\in \A$ with $\mu(A)<\infty$,
\item for every $f\in E$ and every set $A\in \A$ with $\mu(A)<\infty$ $f$ is $\mu$-integrable over $A$,
\item if $g$ is measurable and $f\in E$ such that $\abs*{g(t)}\leq\abs*{f(t)}$ $\mu$-a.\,e. then $g\in E$
and $\norm{g}\leq\norm{f}$.
\end{enumerate}
The standard examples are of course the spaces $L^p(\mu)$ for $1\leq p\leq\infty$.\par
Every K\"othe function space $E$ is a Banach lattice when endowed with the natural order 
$f\leq g$ \ifif $f(t)\leq g(t)$ $\mu$-a.\,e.\par
Recall that a Banach lattice $E$ is said to be order complete ($\sigma$-order complete) if for every net (sequence) in $E$ 
which is order bounded the supremum of said net (sequence) in $E$ exists. A Banach lattice $E$ is called order continuous 
($\sigma$-order continuous) provided that every decreasing net (sequence) in $E$ whose infimum is zero is norm-convergent 
to zero.\par
It is easy to see that a K\"othe function space $E$ is always $\sigma$-order complete and thus by \cite{lin}*{Proposition 3.1.5}
$E$ is order continuous \ifif $E$ is $\sigma$-order continuous \ifif $E$ is order complete and order continuous. Also, reflexivity 
of $E$ implies order continuity, for any $\sigma$-order complete Banach lattice which is not $\sigma$-order continuous contains an 
isomorphic copy of $\ell^{\infty}$ (cf. \cite{lin}*{Proposition 3.1.4}).\par
Let us also mention the following well-known fact that will be needed later.
\begin{lemma}\label{lemma:conv ae}
If $E$ is a K\"othe function space, $(f_n)_{n\in \N}$ a sequence in $E$ and $f\in E$ such that 
$\norm{f_n-f}_E\to 0$ then there is a subsequence of $(f_n)_{n\in \N}$ which converges pointwise 
almost everywhere to $f$.
\end{lemma}
For a K\"othe function space $E$ we denote by $E^{\prime}$ the space of all measurable functions $g:S\rightarrow \R$ (modulo
equality $\mu$-a.\,e.) such that 
\begin{equation*}
\norm{g}_{E^{\prime}}:=\sup\set*{\int_S\abs*{fg}\,\mathrm{d}\mu:f\in B_E}<\infty.
\end{equation*}
Then $(E^{\prime},\norm{\cdot}_{E^{\prime}})$ is again a K\"othe function space, the so called K\"othe dual of $E$. The 
operator $T:E^{\prime}\rightarrow E^*$ defined by 
\begin{equation*}
(Tg)(f)=\int_S fg\,\mathrm{d}\mu \ \ \forall f\in E, \forall g\in E^{\prime}
\end{equation*}
is well-defined, linear and isometric. Moreover, $T$ is onto \ifif $E$ is order continuous (cf.\,\cite{lin}*{p.149}),
thus for order continuous $E$ we have $E^*=E^{\prime}$.\par
We refer the reader to \cite{lindenstrauss} or \cite{lin} for more information on Banach lattices in general and K\"othe 
function spaces in particular.\par
Now recall that if $X$ is a Banach space a function $f:S\rightarrow X$ is called simple if there are finitely many measurable 
sets $A_1,\dots ,A_n\in \A$ such that $\bigcup_{i=1}^{\infty}A_i=S$ and $f$ is constant on each $A_i$. The function $f$
is said to be Bochner-measurable if there exists a sequence $(f_n)_{n\in \N}$ of simple functions such that 
$\lim_{n\to \infty}\norm{f_n(t)-f(t)}=0$ $\mu$-a.\,e. and weakly measurable if $x^*\circ f$ is measurable for every
functional $x^*\in X^*$. According to Pettis' measurability theorem (cf.\,\cite{lin}*{Theorem 3.2.2}) $f$ is Bochner-measurable
\ifif $f$ is weakly measurable and almost everywhere separably valued (i.\,e. there is a separable subspace $Y\ssq X$ such
that $f(t)\in Y$ $\mu$-a.\,e.).\par
For a K\"othe function space $E$ and a Banach space $X$ we denote by $E(X)$ the space of all Bochner-measurable functions
$f:S\rightarrow X$ (modulo equality a.\,e.) such that $\norm{f(\cdot)}\in E$. Endowed with the norm $\norm{f}_{E(X)}=
\norm*{\norm{f(\cdot)}}_E$ $E(X)$ becomes a Banach space, the so called K\"othe Bochner space induced by $E$ and $X$.
The most prominent examples are again the Lebesgue-Bochner spaces $L^p(X)$ for $1\leq p\leq\infty$.\par
Next we recall how the dual of $E(X)$ can be described provided that $E$ is order continuous. A function $F:S\rightarrow X^*$ 
is called weak*-measurable if $F(\cdot)(x)$ is measurable for every $x\in X$. We define an equivalence relation on the set of
all weak*-measurable functions by setting $F\sim G$ \ifif $F(t)(x)=G(t)(x)$ a.\,e. and we write $E^{\prime}(X^*,w^*)$ for the
space of all (equivalence classes of) weak*-measurable functions $F$ such that there is some $g\in E^{\prime}$ with $\norm{F(t)}
\leq g(t)$ a.\,e.\par
A norm on $E^{\prime}(X^*,w^*)$ can be defined by
\begin{equation*}
\norm{[F]}_{E^{\prime}(X^*,w^*)}:=\inf\set*{\norm{g}_{E^{\prime}}:g\in E^{\prime} \ \text{and} \ \norm{F(t)}\leq g(t) \ \text{a.\,e.}}.
\end{equation*}
Then the following deep theorem holds.
\begin{theorem}[cf. \cite{bukhvalov}]\label{thm:dual E(X)}
Let $E$ be an order continuous K\"othe function space over the complete, $\sigma$-finite measure space $(S,\A,\mu)$ and 
let $X$ be a Banach space. Then the map $V:E^{\prime}(X^*,w^*)\rightarrow E(X)^*$ defined by
\begin{equation*}
V([F])(f):=\int_S F(t)(f(t))\,\mathrm{d}\mu(t) \ \ \forall f\in E(X), \forall [F]\in E^{\prime}(X^*,w^*)
\end{equation*}
is an isometric isomorphism and moreover every equivalence class $L$ in $E^{\prime}(X^*,w^*)$ has a representative $F$
such that $\norm{F(\cdot)}\in E^{\prime}$ and $\norm{L}_{E^{\prime}(X^*,w^*)}=\norm*{\norm{F(\cdot)}}_{E^{\prime}}$.
\end{theorem}
Sirotkin proved in \cite{sirotkin} that for $1<p<\infty$ the Lebesgue-Bochner space $L^p(X)$ is \ac{acs} resp. \ac{luacs}
resp. \ac{uacs} whenever $X$ has the respective property. In the next section we will study the more general case of 
K\"othe-Bochner spaces.

\section{Results and proofs}\label{sec:results}
We begin with the \ac{acs} spaces, for which we have the following result.
\begin{proposition}\label{prop:acs}
If $E$ is an order continuous \ac{acs} K\"othe function space and $X$ is an \ac{acs} Banach space,
then $E(X)$ is \ac{acs} as well.
\end{proposition}

\begin{Proof}
The proof is similar to that of \cite{hardtke}*{Proposition 3.3}. First we fix two elements $f,g\in S_{E(X)}$
such that $\norm{f+g}_{E(X)}=2$ and a functional $l\in S_{E(X)^*}$ with $l(f)=1$.\par
Since $E$ is order continuous, by Theorem \ref{thm:dual E(X)} $l$ can be represented via an element 
$[F]\in E^{\prime}(X^*,w^*)$ such that $\norm{F(\cdot)}\in E^{\prime}$ and $\norm*{\norm{F(\cdot)}}_{E^{\prime}}=
\norm{[F]}_{E^{\prime}(X^*,w^*)}=\norm{l}=1$. It follows that
\begin{align*}
&1=l(f)=\int_S F(t)(f(t))\,\mathrm{d}\mu(t)\leq\int_S \norm{F(t)}\norm{f(t)}\,\mathrm{d}\mu(t) \\
&\leq\norm*{\norm{F(\cdot)}}_{E^{\prime}}\norm*{\norm{f(\cdot)}}_E=\norm{l}\norm{f}_{E(X)}=1
\end{align*}
and hence
\begin{equation}\label{eq:3.1}
\int_S \norm{F(t)}\norm{f(t)}\,\mathrm{d}\mu(t)=1
\end{equation}
and
\begin{equation}\label{eq:3.2}
F(t)(f(t))=\norm{F(t)}\norm{f(t)} \ \ \text{a.\,e.}
\end{equation}
We also have
\begin{align*}
&2=\norm{f+g}_{E(X)}=\norm*{\norm{f(\cdot)+g(\cdot)}}_E\leq\norm*{\norm{f(\cdot)}+\norm{g(\cdot)}}_E \\
&\leq\norm{f}_{E(X)}+\norm{g}_{E(X)}=2
\end{align*}
and thus
\begin{equation}\label{eq:3.3}
\norm*{\norm{f(\cdot)}+\norm{g(\cdot)}}_E=2.
\end{equation}
Since $E$ is \ac{acs} it follows from \eqref{eq:3.1} and \eqref{eq:3.3} that
\begin{equation}\label{eq:3.4}
\int_S \norm{F(t)}\norm{g(t)}\,\mathrm{d}\mu(t)=1.
\end{equation}
In a similar way as we have obtained \eqref{eq:3.3} we can also show
\begin{equation}\label{eq:3.5}
\norm*{\norm{f(\cdot)+g(\cdot)}+\norm{f(\cdot)}+\norm{g(\cdot)}}_E=4.
\end{equation}
Because $E$ is \ac{acs} this together with \eqref{eq:3.1}, \eqref{eq:3.3} and \eqref{eq:3.4} implies
\begin{equation}\label{eq:3.6}
\int_S \norm{F(t)}\norm{f(t)+g(t)}\,\mathrm{d}\mu(t)=2.
\end{equation}
From \eqref{eq:3.1}, \eqref{eq:3.4} and \eqref{eq:3.6} we get
\begin{equation}\label{eq:3.7}
\norm{F(t)}\paren*{\norm{f(t)}+\norm{g(t)}-\norm{f(t)+g(t)}}=0 \ \ \text{a.\,e.}
\end{equation}
Now we will show that
\begin{equation}\label{eq:3.8}
F(t)(g(t))=\norm{F(t)}\norm{g(t)} \ \ \text{a.\,e.}
\end{equation}
To this end, let us denote by $N_1$ resp. $N_2$ the null sets on which the equality from \eqref{eq:3.2} 
resp. \eqref{eq:3.7} does not hold. Let $N=N_1\cup N_2$.\par
Put $B=\set*{t\in S\sm N:F(t)\neq 0 \ \text{and} \ g(t)\neq 0}$ and $C=\set*{t\in B:f(t)=0}$. We claim 
that $C$ is a null set.\par
To see this, define $h:S\rightarrow \R$ by $h(t)=\norm{F(t)}$ for $t\in S\sm C$ and $h(t)=0$ for $t\in C$.
Then $h$ is measurable and since $h(t)\leq\norm{F(t)}$ for all $t\in S$ we have $h\in E^{\prime}$ with
$\norm{h}_{E^{\prime}}\leq 1$. We also have $h(t)\norm{f(t)}=\norm{F(t)}\norm{f(t)}$ for every $t\in S$
and hence by \eqref{eq:3.1} 
\begin{equation*}
\int_S h(t)\norm{f(t)}\,\mathrm{d}\mu(t)=1,
\end{equation*}
which also implies $\norm{h}_{E^{\prime}}=1$. Together with \eqref{eq:3.3} we now get
\begin{equation*}
\int_S h(t)\norm{g(t)}\,\mathrm{d}\mu(t)=1,
\end{equation*}
since $E$ is \ac{acs}. Taking into account \eqref{eq:3.4} we arrive at
\begin{equation*}
\int_S \paren*{\norm{F(t)}-h(t)}\norm{g(t)}\,\mathrm{d}\mu(t)=0.
\end{equation*}
Hence $\paren*{\norm{F(t)}-h(t)}\norm{g(t)}=0$ a.\,e. and thus $C$ must be a null set.\par
Now if $t\in (S\sm C)\cap B$ then $F(t)\neq 0$, $f(t)\neq 0$ and $g(t)\neq 0$ and 
$\norm{F(t)}\norm{f(t)}=F(t)(f(t))$ as well as
\begin{equation*}
\norm{f(t)+g(t)}=\norm{f(t)}+\norm{g(t)}.
\end{equation*}
By \cite{abramovich}*{Lemma 2.1} this implies
\begin{equation*}
\norm*{\frac{f(t)}{\norm{f(t)}}+\frac{g(t)}{\norm{g(t)}}}=2.
\end{equation*}
Since $X$ is \ac{acs} it follows that $\norm{F(t)}\norm{g(t)}=F(t)(g(t))$.\par
So $M:=N\cup C$ is a null set with $\norm{F(t)}\norm{g(t)}=F(t)(g(t))$ for every $t\in S\sm M$
and \eqref{eq:3.8} is proved.\par
Now combining \eqref{eq:3.4} and \eqref{eq:3.8} we obtain
\begin{equation*}
l(g)=\int_S F(t)(g(t))\,\mathrm{d}\mu(t)=1,
\end{equation*}
which finishes the proof.
\end{Proof}

Before we turn to the case of \ac{luacs} spaces, let us recall Egorov's theorem (cf. \cite{werner2}*{Satz IV.6.7}), 
which states that for any finite measure space $(S,\A,\mu)$ and every sequence $(f_n)_{n\in \N}$ of measurable functions 
on $S$ which converges to zero pointwise $\mu$-a.\,e. and each $\eps>0$ there is a set $A\in \A$ with 
$\mu(S\sm A)\leq\eps$ such that $(f_n)_{n\in \N}$ is uniformly convergent to zero on $A$.\par
Now we are ready to prove the following theorem.
\begin{theorem}\label{thm:luacs}
Let $E$ be an order continuous K\"othe function space over the complete $\sigma$-finite measure space 
$(S,\A,\mu)$ and $X$ an \ac{luacs} Banach space. If  
\begin{enumerate}[\upshape(a)]
\item $E$ is \ac{WLUR} or
\item $E$ is $\text{luacs}^+$ and $E^{\prime}$ is also order continuous
\end{enumerate}
then $E(X)$ is also \ac{luacs}.
\end{theorem}

\begin{Proof}
Suppose that we are given a sequence $(f_n)_{n\in \N}$ in $S_{E(X)}$ and an element $f\in S_{E(X)}$
such that $\norm{f_n+f}_{E(X)}\to 2$ as well as a functional $l\in S_{E(X)^*}$ such that $l(f_n)\to 1$.
As before, we can represent $l$ by an element $[F]\in E^{\prime}(X^*,w^*)$. We then have
\begin{equation*}
l(f_n)=\int_S F(t)(f_n(t))\,\mathrm{d}\mu(t)\leq\int_S \norm{F(t)}\norm{f_n(t)}\,\mathrm{d}\mu(t)\leq 1
\end{equation*}
and hence
\begin{equation}\label{eq:3.9}
\lim_{n\to \infty}\int_S \norm{F(t)}\norm{f_n(t)}\,\mathrm{d}\mu(t)=1.
\end{equation}
By passing to a subsequence we may also assume that
\begin{equation}\label{eq:3.10}
\lim_{n\to \infty}\paren*{\norm{F(t)}\norm{f_n(t)}-F(t)(f_n(t))}=0 \ \ \text{a.\,e.}
\end{equation}
We further have
\begin{equation*}
\norm{f_n+f}_{E(X)}=\norm*{\norm{f_n(\cdot)+f(\cdot)}}_E\leq\norm*{\norm{f_n(\cdot)}+\norm{f(\cdot)}}_E\leq2
\end{equation*}
and thus
\begin{equation}\label{eq:3.11}
\lim_{n\to \infty}\norm*{\norm{f_n(\cdot)}+\norm{f(\cdot)}}_E=2.
\end{equation}
An analogous argument also shows 
\begin{equation}\label{eq:3.12}
\lim_{n\to \infty}\norm*{\norm{f_n(\cdot)+f(\cdot)}+\norm{f_n(\cdot)}+\norm{f(\cdot)}}_E=4.
\end{equation}
Moreover, the inequality
\begin{align*}
&\norm{f_n+f}_{E(X)}+1\geq\norm*{\norm{f_n(\cdot)+f(\cdot)}+\norm{f_n(\cdot)}}_E \\
&\geq\norm*{\norm{f_n(\cdot)+f(\cdot)}+\norm{f_n(\cdot)}+\norm{f(\cdot)}}_E-1
\end{align*}
holds for every $n\in \N$. It follows that
\begin{equation}\label{eq:3.13}
\lim_{n\to \infty}\norm*{\norm{f_n(\cdot)+f(\cdot)}+\norm{f_n(\cdot)}}_E=3.
\end{equation}
Analogously one can see that
\begin{equation}\label{eq:3.14}
\lim_{n\to \infty}\norm*{\norm{f_n(\cdot)+f(\cdot)}+\norm{f(\cdot)}}_E=3.
\end{equation}
Finally, we have
\begin{align*}
&\norm*{\norm{f_n(\cdot)+f(\cdot)}+\norm{f_n(\cdot)}}_E+3\geq\norm*{\norm{f_n(\cdot)+f(\cdot)}+\norm{f_n(\cdot)}+3\norm{f(\cdot)}}_E \\
&\geq2\norm*{\norm{f_n(\cdot)+f(\cdot)}+\norm{f(\cdot)}}_E,
\end{align*}
consequently
\begin{equation}\label{eq:3.15}
\lim_{n\to \infty}\norm*{\norm{f_n(\cdot)+f(\cdot)}+\norm{f_n(\cdot)}+3\norm{f(\cdot)}}_E=6.
\end{equation}
Since $E$ is in particular \ac{luacs} we get from \eqref{eq:3.9} and \eqref{eq:3.11} that
\begin{equation}\label{eq:3.16}
\int_S \norm{F(t)}\norm{f(t)}\,\mathrm{d}\mu(t)=1.
\end{equation}
Because $E$ is in any case $\mathrm{luacs}^+$ it follows from \eqref{eq:3.13}, \eqref{eq:3.15} and \eqref{eq:3.16} that
\begin{equation*}
\lim_{n\to \infty}\int_S \norm{F(t)}\paren*{\norm{f_n(t)}+\norm{f_n(t)+f(t)}}\,\mathrm{d}\mu(t)=3.
\end{equation*}
and thus
\begin{equation*}
\lim_{n\to \infty}\int_S \norm{F(t)}\paren*{\norm{f_n(t)}+\norm{f(t)}-\norm{f_n(t)+f(t)}}\,\mathrm{d}\mu(t)=0.
\end{equation*}
So by passing to a further subsequence we may assume 
\begin{equation}\label{eq:3.17}
\lim_{n\to \infty}\norm{F(t)}\paren*{\norm{f_n(t)}+\norm{f(t)}-\norm{f_n(t)+f(t)}}=0 \ \ \text{a.\,e.}
\end{equation}
Next we will show that 
\begin{equation}\label{eq:3.18}
F(t)(f(t))=\norm{F(t)}\norm{f(t)} \ \ \text{a.\,e.}
\end{equation}
Since $(S,\A,\mu)$ is $\sigma$-finite there is an increasing sequence $(A_m)_{m\in \N}$ in $\A$ such that $\mu(A_m)<\infty$
for every $m\in \N$ and $\bigcup_{m=1}^{\infty}A_m=S$.\par
Denote by $N_1$ resp. $N_2$ the null sets on which the convergence statement from \eqref{eq:3.10} resp. \eqref{eq:3.17} does 
not hold and let $N=N_1\cup N_2$. Put $B=\set*{t\in S\sm N:F(t)\neq 0 \ \mathrm{and} \ f(t)\neq 0}$ and 
$C=\set*{t\in B:\norm{f_n(t)}\to 0}$. We shall see that $C$ is a null set.\par
First we define for every $m\in \N$ a function $a_m:S\rightarrow \R$ by setting $a_m(t)=\norm{F(t)}$ for $t\in S\sm (C\cap A_m)$ 
and $a_m(t)=0$ for $t\in C\cap A_m$. Note that each $a_m$ is measurable and since $\abs*{a_m(t)}\leq \norm{F(t)}$ for every $t\in S$ 
we have $a_m\in B_{E^{\prime}}$.\par
We have $\lim_{k\to \infty}\norm{F(t)}\norm{f_k(t)}\chi_{C\cap{A_m}}(t)=0$ for every $t\in S$ and every $m\in \N$, so by Egorov's theorem
we can find for every $m\in \N$ an increasing sequence $(B_{n,m})_{n\in \N}$ in $\A|_{A_m}$ with $\mu(A_m\sm B_{n,m})\leq 1/n$ and 
such that $(\norm{F(\cdot)}\norm{f_k(\cdot)}\chi_{C\cap{A_m}})_{k\in \N}$ converges uniformly to zero on each $B_{n,m}$.\par
It follows that $M_m:=\bigcap_{n=1}^{\infty}A_m\sm B_{n,m}$ is a null set for every $m\in \N$.\par
Let us now first suppose that (b) holds, so $E^{\prime}$ is order continuous. We have
\begin{equation*}
\lim_{n\to \infty}\norm{F(t)}\chi_{C\cap (A_m\sm B_{n,m})}(t)=0 \ \ \forall t\in S\sm M_m
\end{equation*}
and moreover this sequence is decreasing, so the order continuity of $E^{\prime}$ implies
\begin{equation*} 
\lim_{n\to \infty}\norm{\norm{F(\cdot)}\chi_{C\cap (A_m\sm B_{n,m})}}_{E^{\prime}}=0.
\end{equation*}
So if $m\in \N$ and $\eps>0$ are given we can find an index $n\in \N$ such that $\norm{\norm{F(\cdot)}\chi_{C\cap (A_m\sm B_{n,m})}}_{E^{\prime}}\leq\eps$
and then, by uniform convergence, an index $k_0\in \N$ such that $\norm{F(t)}\norm{f_k(t)}\chi_{C\cap B_{n,m}}(t)\leq\eps\mu(A_m)^{-1}$ for every $t\in S$ 
and every $k\geq k_0$.\par
Then we have 
\begin{align*}
&\int_{C\cap A_m}\norm{F(t)}\norm{f_k(t)}\,\mathrm{d}\mu(t) \\
&=\int_{C\cap B_{n,m}}\norm{F(t)}\norm{f_k(t)}\,\mathrm{d}\mu(t)+\int_{C\cap (A_m\sm B_{n,m})}\norm{F(t)}\norm{f_k(t)}\,\mathrm{d}\mu(t) \\
&\leq\int_{C\cap B_{n,m}} \frac{\eps}{\mu(A_m)}\,\mathrm{d}\mu(t)+\norm{\norm{F(\cdot)}\chi_{C\cap (A_m\sm B_{n,m})}}_{E^{\prime}}\leq2\eps
\end{align*}
for each $k\geq k_0$.\par
In conclusion we have 
\begin{equation*}\label{+}
\lim_{k\to \infty}\int_{C\cap A_m}\norm{F(t)}\norm{f_k(t)}\,\mathrm{d}\mu(t)=0 \ \ \forall m\in \N. \tag{+}
\end{equation*}
Now if (a) holds, i.\,e. if $E$ is \ac{WLUR} then by \eqref{eq:3.11} the sequence $(\norm{f_k(\cdot)})_{k\in \N}$ must be
weakly convergent to $\norm{f(\cdot)}$ in $E$ and hence 
\begin{equation*}
\lim_{k\to \infty}\int_{C\cap (A_m\sm B_{n,m})}\norm{F(t)}\norm{f_k(t)}\,\mathrm{d}\mu(t)=\int_{C\cap (A_m\sm B_{n,m})}\norm{F(t)}\norm{f(t)}\,\mathrm{d}\mu(t)
\end{equation*}
for all $n,m\in \N$. Since $(\norm{f(\cdot)}\chi_{C\cap (A_m\sm B_{n,m})})_{n\in \N}$ dereases to zero a.\,e. the 
order continuity of $E$ gives us $\lim_{n\to \infty}\norm{\norm{f(\cdot)}\chi_{C\cap (A_m\sm B_{n,m})}}_E=0$ for every $m\in \N$.\par
A similiar argument as before now easily yields that \eqref{+} also holds in case (a). But \eqref{+} is nothing else than
\begin{equation*}
\lim_{n\to \infty}\int_S \paren*{\norm{F(t)}-a_m(t)}\norm{f_n(t)}\,\mathrm{d}\mu(t)=0 \ \ \forall m\in \N.
\end{equation*}
Combinig this with \eqref{eq:3.9} leaves us with
\begin{equation*}
\lim_{n\to \infty}\int_S a_m(t)\norm{f_n(t)}\,\mathrm{d}\mu(t)=1 \ \ \forall m\in \N.
\end{equation*}
Since $E$ is \ac{luacs} and because of \eqref{eq:3.11} it follows that
\begin{equation*}
\int_S a_m(t)\norm{f(t)}\,\mathrm{d}\mu(t)=1 \ \ \forall m\in \N.
\end{equation*}
Taking into account \eqref{eq:3.16} we get
\begin{equation*}
\int_S \paren*{\norm{F(t)}-a_m(t)}\norm{f(t)}\,\mathrm{d}\mu(t)=0 \ \ \forall m\in \N
\end{equation*}
and hence for every $m\in \N$ we have $\paren*{\norm{F(t)}-a_m(t)}\norm{f(t)}=0$ a.\,e.
Consequently, $C\cap A_m$ is a null set for every $m$ and thus $C=\bigcup_{m=1}^{\infty}C\cap A_m$ is also a null set.\par
Now suppose that $t\in (S\sm C)\cap B$. Then we have $F(t)\neq 0$, $f(t)\neq 0$ and $\norm{f_n(t)}\not\to 0$, as well as
$\norm{F(t)}\norm{f_n(t)}-F(t)(f_n(t))\to 0$ and
\begin{equation*}
\lim_{n\to \infty}\paren*{\norm{f_n(t)}+\norm{f(t)}-\norm{f_n(t)+f(t)}}=0.
\end{equation*}
By passing to a subsequence we may assume that $(\norm{f_n(t)})_{n\in \N}$ is bounded away from zero. Then it follows
from Lemma \ref{aux lemma} that
\begin{equation*}
\lim_{n\to \infty}\norm*{\frac{f_n(t)}{\norm{f_n(t)}}+\frac{f(t)}{\norm{f(t)}}}=2.
\end{equation*}
Also, we have
\begin{equation*}
\lim_{n\to \infty}\frac{F(t)}{\norm{F(t)}}\paren*{\frac{f_n(t)}{\norm{f_n(t)}}}=1.
\end{equation*}
Since $X$ is \ac{luacs} we can conclude that $F(t)(f(t))=\norm{F(t)}\norm{f(t)}$.\par
So $M:=N\cup C$ is a null set with $F(t)(f(t))=\norm{F(t)}\norm{f(t)}$ for every $t\in S\sm M$ and
\eqref{eq:3.18} is proved.\par
From \eqref{eq:3.16} and \eqref{eq:3.18} it follows that
\begin{equation*}
l(f)=\int_S F(t)(f(t))\,\mathrm{d}\mu(t)=1
\end{equation*}
and we are done.
\end{Proof}

Recall that a subset $A\ssq L^1(\mu)$ is said to be equi-integrable if for every $\eps>0$ there is some 
$\delta>0$ such that 
\begin{equation*}
B\in \A \ \mathrm{with} \ \mu(B)\leq\delta \ \Rightarrow \ \abs*{\int_B f\,\mathrm{d}\mu}\leq\eps \ \ \forall f\in A.
\end{equation*}
It is well-known that for a finite measure $\mu$ a bounded subset $A\ssq L^1(\mu)$ is relatively weakly compact 
in $L^1(\mu)$ if and only if $A$ is equi-integrable (see for instance \cite{werner1}*{Satz VIII.6.9}). One ingredient 
for the usual proof of this fact is the following lemma (see \cite{werner1}*{Lemma VIII.6.7}), which we will also need in the sequel.
\begin{lemma}\label{lemma:equi-int}
For a finite measure space $(S,\A,\mu)$, a sequence $(f_n)_{n\in \N}$ in $L^1(\mu)$ is equi-integrable 
whenever the sequence $(\int_B f_n\,\mathrm{d}\mu)_{n\in \N}$ is convergent for each $B\in \A$.
\end{lemma}
We will also need Vitali's Lemma, which reads as follows (see for example \cite{lin}*{Lemma 3.1.13}) for an even more general version).
\begin{lemma}\label{lemma:vitali}
Let $(S,\A,\mu)$ be a finite measure space and let $(f_n)_{n\in \N}$ be a sequence in $L^1(\mu)$ 
such that $\set*{\abs*{f_n}:n\in \N}$ is equi-integrable. Let $f$ be a measurable function on $S$ 
such that $f_n(t)\to f(t)$ $\mu$-a.\,e. Then $f\in L^1(\mu)$ and $\norm{f_n-f}_1\to 0$.
\end{lemma}
Finally, let us recall that a Banach space $X$ is said to have the Kadets-Klee property (also known as 
property $(H)$) if for every sequence $(x_n)_{n\in \N}$ in $X$ and each $x\in X$ the implication 
\begin{equation*}
x_n\xrightarrow{\sigma} x \ \mathrm{and} \ \norm{x_n}\to \norm{x} \ \Rightarrow \ \norm{x_n-x}\to 0
\end{equation*}
holds. For example, every \ac{LUR} space and every dual of a reflexive, \ac{FS} space has the 
Kadets-Klee property.\par
It is known that every Banach lattice with the Kadets-Klee property is order 
continuous (cf.\,\cite{lindenstrauss}*{p.28}). With this in mind we can prove 
the following result concerning $\mathrm{luacs}^+$ spaces.
\begin{theorem}\label{thm:luacs+}
If the measure $\mu$ is finite and $E$ is \ac{LUR}, then $E(X)$ is a $\text{luacs}^+$ space 
whenever $X$ is $\text{luacs}^+$. If in addition $E^{\prime}$ is order continuous then the 
assertion also holds if $\mu$ is merely $\sigma$-finite.
\end{theorem}

\begin{Proof}
By the previous theorem, $E(X)$ is \ac{luacs} so we only have to show the implication ``$\Leftarrow$''
in Definition \ref{def:luacs+sluacs+} (i). To this end, let $(f_n)_{n\in \N}$ be a sequence in $S_{E(X)}$
and $f\in S_{E(X)}$ such that $\norm{f_n+f}_{E(X)}\to 2$ and let $l\in S_{E(X)^*}$ such that $l(f)=1$.
It will be enough to show that a subsequence of $(l(f_n))_{n\in \N}$ converges to one.\par
Since $E$ is order continuous we can as before represent $l$ by some $[F]\in E^{\prime}(X^*,w^*)$ 
and conclude
\begin{equation}\label{eq:3.19}
\int_S\norm{F(t)}\norm{f(t)}\,\mathrm{d}\mu(t)=1
\end{equation}
and
\begin{equation}\label{eq:3.20}
\norm{F(t)}\norm{f(t)}=F(t)(f(t)) \ \ \text{a.\,e.}
\end{equation}
Also, just as we have done in the previous proof, we find that
\begin{align}
&\lim_{n\to \infty}\norm*{\norm{f_n(\cdot)}+\norm{f(\cdot)}}_E=2,\label{eq:3.21} \\
&\lim_{n\to \infty}\norm*{\norm{f_n(\cdot)+f(\cdot)}+\norm{f_n(\cdot)}+\norm{f(\cdot)}}_E=4,\label{eq:3.22} \\
&\lim_{n\to \infty}\norm*{\norm{f_n(\cdot)+f(\cdot)}+\norm{f_n(\cdot)}}_E=3,\label{eq.3.23} \\
&\lim_{n\to \infty}\norm*{\norm{f_n(\cdot)+f(\cdot)}+\norm{f(\cdot)}}_E=3,\label{eq:3.24} \\
&\lim_{n\to \infty}\norm*{\norm{f_n(\cdot)+f(\cdot)}+\norm{f_n(\cdot)}+3\norm{f(\cdot)}}_E=6.\label{eq:3.25}
\end{align}
Since $E$ is \ac{LUR} it follows that
\begin{align}
&\lim_{n\to \infty}\norm*{\norm{f_n(\cdot)}-\norm{f(\cdot)}}_E=0,\label{eq:3.26} \\
&\lim_{n\to \infty}\norm*{\norm{f_n(\cdot)+f(\cdot)}+\norm{f_n(\cdot)}-3\norm{f(\cdot)}}_E=0.\label{eq:3.27}
\end{align}
Hence by passing to a subsequence we may assume that (cf.\,Lemma \ref{lemma:conv ae})
\begin{align}
&\lim_{n\to \infty}\norm{f_n(t)}=\norm{f(t)} \ \ \text{a.\,e.} \ \ \text{and}\label{eq:3.28} \\
&\lim_{n\to \infty}\norm{f_n(t)+f(t)}=2\norm{f(t)} \ \ \text{a.\,e.}\label{eq:3.29}
\end{align}
By \eqref{eq:3.26} and \eqref{eq:3.19} we also have
\begin{equation}\label{eq:3.30}
\lim_{n\to \infty}\int_S\norm{F(t)}\norm{f_n(t)}\,\mathrm{d}\mu(t)=1.
\end{equation}
Since $X$ is $\mathrm{luacs}^+$ it follows from \eqref{eq:3.20}, \eqref{eq:3.28} and \eqref{eq:3.29} that
\begin{equation}\label{eq:3.31}
\lim_{n\to \infty}\paren*{\norm{F(t)}\norm{f_n(t)}-F(t)(f_n(t))}=0 \ \ \text{a.\,e.}
\end{equation}
From \eqref{eq:3.26} we also get 
\begin{equation*}
\lim_{n\to \infty}\int_A\norm{F(t)}\norm{f_n(t)}\,\mathrm{d}\mu(t)=\int_A\norm{F(t)}\norm{f(t)}\,\mathrm{d}\mu(t) \ \ \forall A\in \A.
\end{equation*}
Thus by Lemma \ref{lemma:equi-int} the sequence $(\norm{F(\cdot)}\norm{f_n(\cdot)}\chi_B)_{n\in \N}$ and hence also the sequence
$((\norm{F(\cdot)}\norm{f_n(\cdot)}-F(\cdot)(f_n(\cdot)))\chi_B)_{n\in \N}$ is equi-integrable with respect to $(B,\A|_B,\mu_{\A|_B})$ for every 
$B\in \A$ with $\mu(B)<\infty$. This combined with Vitali's Lemma and \eqref{eq:3.31} implies
\begin{equation*}
\lim_{n\to \infty}\int_B\paren*{\norm{F(t)}\norm{f_n(t)}-F(t)(f_n(t))}\,\mathrm{d}\mu(t)=0 \ \ \forall B\in \A \ \mathrm{with} \ \mu(B)<\infty.
\end{equation*}
So if $\mu(S)<\infty$ we immediately get 
\begin{equation*}
l(f_n)=\int_S F(t)(f_n(t))\,\mathrm{d}\mu(t)\to 1,
\end{equation*}
because of \eqref{eq:3.30}.\par
If $\mu$ is merely $\sigma$-finite but $E^{\prime}$ is order continuous, we can fix an increasing sequence $(A_m)_{m\in \N}$ in $\A$
such that $\bigcup_{m=1}^{\infty}A_m=S$ and $\mu(A_m)<\infty$ for every $m\in \N$. Then the sequence $(\norm{F(\cdot)}\chi_{S\sm A_m})_{m\in \N}$
decreases pointwise to zero and, by the order continuity of $E^{\prime}$, we can conclude that $\norm{\norm{F(\cdot)}\chi_{S\sm A_m}}_{E^{\prime}}\to 0$.\par
Thus given any $\eps>0$ we find an $m_0\in \N$ such that $\norm{\norm{F(\cdot)}\chi_{S\sm A_{m_0}}}_{E^{\prime}}\leq\eps/3$. Since $\mu(A_{m_0})<\infty$
there exists $N\in \N$ such that
\begin{equation*}
\int_{A_{m_0}}\paren*{\norm{F(t)}\norm{f_n(t)}-F(t)(f_n(t))}\,\mathrm{d}\mu(t)\leq\frac{\eps}{3} \ \ \forall n\geq N.
\end{equation*}
It follows that for every $n\geq N$
\begin{align*}
&\int_S\paren*{\norm{F(t)}\norm{f_n(t)}-F(t)(f_n(t))}\,\mathrm{d}\mu(t)\leq\frac{\eps}{3}+2\int_{S\sm A_{m_0}}\norm{F(t)}\norm{f_n(t)}\,\mathrm{d}\mu(t) \\
&\leq\frac{\eps}{3}+2\norm{\norm{F(\cdot)}\chi_{S\sm A_{m_0}}}_{E^{\prime}}\leq\eps.
\end{align*}
So we have
\begin{equation*}
\lim_{n\to \infty}\int_S\paren*{\norm{F(t)}\norm{f_n(t)}-F(t)(f_n(t))}\,\mathrm{d}\mu(t)=0
\end{equation*}
and because of \eqref{eq:3.30} it follows as before that
\begin{equation*}
l(f_n)=\int_S F(t)(f_n(t))\,\mathrm{d}\mu(t)\to 1,
\end{equation*}
finishing the proof.
\end{Proof}

Now we turn to the \ac{sluacs} spaces. An easy normalisation argument shows that a Banach space $X$ is 
\ac{sluacs} \ifif for every $x\in S_X$, every sequence $(x_n^*)_{n\in \N}$ in $S_{X^*}$ and all sequences 
$(x_n)_{n\in \N}$ in $X$ with $\norm{x_n+x}\to 2$, $\norm{x_n}\to 1$ and $x_n^*(x_n)\to 1$ we have 
$x_n^*(x)\to 1$. In view of this characterisation, $X$ is \ac{sluacs} \ifif for every $x\in S_X$ and
every $0<\eps\leq 2$ the number
\begin{equation*}
\beta_X(x,\eps):=\inf\set*{\max\set*{1-\norm*{\frac{x+y}{2}},\abs*{\norm{y}-1},\abs*{x^*(y)-1}}:(y,x^*)\in V_{x,\eps}}
\end{equation*}
is strictly positive, where
\begin{equation*}
V_{x,\eps}:=\set*{(y,x^*)\in X\times S_{X^*}:x^*(y-x)\geq\eps}.
\end{equation*}
Next we will prove an easy Lemma on the continuity of $\beta_X$.
\begin{lemma}\label{lemma:cont beta}
For all $0<\eps,\tilde{\eps},\leq 2$ and all $x,\tilde{x}\in S_X$ we have
\begin{equation*}
\abs*{\beta_X(x,\eps)-\beta_X(\tilde{x},\tilde{\eps})}\leq\norm{x-\tilde{x}}+\abs{\eps-\tilde{\eps}},
\end{equation*}
i.\,e. $\beta_X$ is $1$-Lipschitz continuous with respect to the norm of $X\oplus_1 \R$.
\end{lemma}

\begin{Proof}
First we fix $0<\eps\leq2$ and $x,\tilde{x}\in S_X$. Put $\delta=\norm{x-\tilde{x}}$ and take 
$y\in X$, $x^*\in S_{X^*}$ such that $x^*(y-x)\geq\eps$. It follows that $x^*(y-\tilde{x})\geq
\eps-\delta$.\par
Now let $0<\tau<1$ be arbitrary. We can find $z\in S_X$ with $x^*(z)\geq 1-\tau$. Define 
$\tilde{y}=y+\delta(1-\tau)^{-1}z$. Then 
\begin{equation*}
x^*(\tilde{y}-\tilde{x})=\frac{\delta}{1-\tau}x^*(z)+x^*(y-\tilde{x})\geq\delta+x^*(y-\tilde{x})=\eps
\end{equation*}
and hence
\begin{equation*}
\max\set*{1-\norm*{\frac{\tilde{x}+\tilde{y}}{2}},\abs*{\norm{\tilde{y}}-1},\abs*{x^*(\tilde{y})-1}}\geq\beta_X(\tilde{x},\eps).
\end{equation*}
But we have $\abs*{\norm{\tilde{y}}-\norm{y}}\leq\norm{y-\tilde{y}}=\delta(1-\tau)^{-1}$ and 
$\abs*{x^*(\tilde{y})-x^*(y)}\leq\norm{y-\tilde{y}}=\delta(1-\tau)^{-1}$ as well as
\begin{equation*}
\abs*{\norm*{\frac{x+y}{2}}-\norm*{\frac{\tilde{x}+\tilde{y}}{2}}}\leq\frac{1}{2}(\norm{x-\tilde{x}}+\norm{y-\tilde{y}})
=\frac{1}{2}\paren*{\delta+\frac{\delta}{1-\tau}}\leq\frac{\delta}{1-\tau}.
\end{equation*}
Thus we get 
\begin{equation*}
\max\set*{1-\norm*{\frac{x+y}{2}},\abs*{\norm{y}-1},\abs*{x^*(y)-1}}\geq\beta_X(\tilde{x},\eps)-\frac{\delta}{1-\tau}
\end{equation*}
and since $0<\tau<1$ was arbitrary it follows that
\begin{equation*}
\max\set*{1-\norm*{\frac{x+y}{2}},\abs*{\norm{y}-1},\abs*{x^*(y)-1}}\geq\beta_X(\tilde{x},\eps)-\delta.
\end{equation*}
Again, since $(y,x^*)\in V_{x,\eps}$ was arbitray we can conclude that
\begin{equation*}
\beta_X(\tilde{x},\eps)-\beta_X(x,\eps)\leq\delta=\norm{x-\tilde{x}}
\end{equation*}
and by symmetry it folows that
\begin{equation*}
\abs{\beta_X(\tilde{x},\eps)-\beta_X(x,\eps)}\leq\norm{x-\tilde{x}}.
\end{equation*}
Analogously one can prove that
\begin{equation*}
\abs{\beta_X(x,\tilde{\eps})-\beta_X(x,\eps)}\leq\abs{\eps-\tilde{\eps}}
\end{equation*}
for all $x\in S_X$ and all $0<\eps,\tilde{\eps},\leq 2$. An application of the triangle inequality
then yields the result.
\end{Proof}

In the paper \cite{kaminska} A. Kami\'nska and B. Turett proved various theorems concerning different rotundity properties
of K\"othe-Bochner spaces. For example, by \cite{kaminska}*{Theorem 5} if $E$ has the so called Fatou property and is \ac{LUR} 
then $E(X)$ is \ac{LUR} whenever $X$ is \ac{LUR}. We will adopt the technique of proof from \cite{kaminska}*{Theorem 5} to show 
the  following result.
\begin{theorem}\label{thm:sluacs}
If $E$ is \ac{LUR} and $X$ is \ac{sluacs} then $E(X)$ is also \ac{sluacs}.
\end{theorem}

\begin{Proof}
Since $E$ is \ac{LUR} it is order continuous.\par
Let $0<\eps\leq2$ and $f\in S_{E(X)}$ be arbitrary and let
\begin{equation*}
A_n:=\set*{t\in S: f(t)\neq 0 \ \mathrm{and} \ \beta_X\paren*{\frac{f(t)}{\norm{f(t)}},\frac{\eps}{8}}\geq\frac{1}{n}}
\end{equation*}
for every $n\in \N$. Since by the previous lemma $\beta_X(\cdot,\eps/8)$ is continuous it follows that the sets $A_n$ 
are measurable. Also, the sequence $(A_n)_{n\in \N}$ is increasing and because $X$ is \ac{sluacs} we have 
$\bigcup_{n=1}^{\infty}A_n=\set*{t\in S:f(t)\neq 0}$, hence $(\norm{f(\cdot)}\chi_{S\sm A_n})_{n\in \N}$ decreases pointwise
to zero. The order continuity of $E$ implies $\norm{\norm{f(\cdot)}\chi_{S\sm A_n}}_E\to 0$ and thus we can find $n_0\in \N$
with
\begin{equation}\label{eq:3.32}
\norm{\norm{f(\cdot)}\chi_{S\sm A_{n_0}}}_E\leq\frac{\eps}{64}.
\end{equation}
Now let us take $g\in S_{E(X)}$ and $l\in S_{E(X)^*}$ with $l(g)=1$ and $l(f)\leq 1-\eps$. Let $l$ be represented 
by $[F]\in E^{\prime}(X^*,w^*)$. As in the proof of \ref{prop:acs} we can conclude
\begin{equation}\label{eq:3.33}
\int_S\norm{F(t)}\norm{g(t)}\,\mathrm{d}\mu(t)=1
\end{equation}
and
\begin{equation}\label{eq:3.34}
\norm{F(t)}\norm{g(t)}=F(t)(g(t)) \ \ \text{a.\,e.}
\end{equation}
Next we define 
\begin{align*}
&C:=\set*{t\in S:F(t)\neq 0} \ \ \mathrm{and}  \\
&B:=\set*{t\in C:F(t)(g(t)-f(t))\geq\frac{\eps}{4}\norm{F(t)}\max\set*{\norm{f(t)},\norm{g(t)}}}.
\end{align*}
Then $B$ is measurable and 
\begin{align*}
&\int_{S\sm B}F(t)(g(t)-f(t))\,\mathrm{d}\mu(t)\leq\frac{\eps}{4}\int_{S\sm B}\norm{F(t)}\max\set*{\norm{f(t)},\norm{g(t)}}\,\mathrm{d}\mu(t) \\
&\leq\frac{\eps}{4}\int_{S\sm B}\norm{F(t)}(\norm{f(t)}+\norm{g(t)})\,\mathrm{d}\mu(t)\leq\frac{\eps}{4}2=\frac{\eps}{2}.
\end{align*}
Since $l(g-f)\geq\eps$ it follows that
\begin{equation}\label{eq:3.35}
\int_BF(t)(g(t)-f(t))\,\mathrm{d}\mu(t)\geq\frac{\eps}{2}.
\end{equation}
Let us fix $0<\eta<\min\set*{\eps/16,1/2n_0}$ such that 
\begin{equation}\label{eq:3.36}
\frac{\eta}{1-\eta}<\frac{2}{n_0}.
\end{equation}
Now consider the sets 
\begin{align*}
&B_1:=\set*{t\in B:\norm{g(t)}<(1-\eta)\norm{f(t)}}, \\
&B_2:=\set*{t\in B:(1-\eta)\norm{f(t)}\leq\norm{g(t)}\leq\norm{f(t)}}, \\
&B_3:=\set*{t\in B:(1-\eta)\norm{g(t)}\leq\norm{f(t)}<\norm{g(t)}}, \\
&B_4:=\set*{t\in B:(1-\eta)\norm{g(t)}>\norm{f(t)}}.
\end{align*}
Then $B_1,\dots,B_4$ are measurable, pairwise disjoint and $\bigcup_{i=1}^{4}B_i=B$. Thus by \eqref{eq:3.35} there exists some
$i\in \set*{1,\dots,4}$ such that
\begin{equation*}
\int_{B_i}F(t)(g(t)-f(t))\,\mathrm{d}\mu(t)\geq\frac{\eps}{8}.
\end{equation*}
If $i=1$ then, since $\norm{g(t)}\leq\norm{f(t)}$ for $t\in B_1$, it follows that 
\begin{equation*}
\int_{B_1}\norm{F(t)}\norm{f(t)}\,\mathrm{d}\mu(t)\geq\frac{\eps}{16}
\end{equation*}
and again by definition of $B_1$ we obtain
\begin{align*}
&\norm*{\norm{g(\cdot)}-\norm{f(\cdot)}}_E=\norm{\abs{\norm{g(\cdot)}-\norm{f(\cdot)}}}_E \\
&\geq\int_{B_1}\norm{F(t)}(\norm{f(t)}-\norm{g(t)})\,\mathrm{d}\mu(t)\geq\eta\int_{B_1}\norm{F(t)}\norm{f(t)}\,\mathrm{d}\mu(t)\geq\eta\frac{\eps}{16}
\end{align*}
and hence
\begin{equation*}
\norm*{\frac{f+g}{2}}_{E(X)}\leq\norm*{\frac{\norm{f(\cdot)}+\norm{g(\cdot)}}{2}}_E\leq 1-\delta_E\paren*{\norm{f(\cdot)},\eta\frac{\eps}{16}}.
\end{equation*}
In the case $i=4$ one can obtain the same statement by an analogous argument. To treat the remaining cases we need
some preliminary considerations.\par
Let us denote by $N$ the null set on which the equality from \eqref{eq:3.34} does not hold and suppose that $t\in B_2\cap A_{n_0}\cap(S\sm N)$.
Then in particular $t\in B$ and $\norm{f(t)}\geq\norm{g(t)}$ and hence 
\begin{equation*}
\frac{F(t)}{\norm{F(t)}}\paren*{\frac{g(t)}{\norm{f(t)}}-\frac{f(t)}{\norm{f(t)}}}\geq\frac{\eps}{4}.
\end{equation*}
Moreover, by the definitions of $B_2$ and $A_{n_0}$ and the choice of $\eta$ we have
\begin{align*}
&\abs*{\norm*{\frac{g(t)}{\norm{f(t)}}}-1}=\abs*{\frac{\norm{g(t)}}{\norm{f(t)}}-1}\leq\eta<\frac{1}{n_0} \\
&\leq\beta_X\paren*{\frac{f(t)}{\norm{f(t)}},\frac{\eps}{8}}\leq\beta_X\paren*{\frac{f(t)}{\norm{f(t)}},\frac{\eps}{4}}.
\end{align*}
Since $t\in (S\sm N)$ we also have
\begin{equation*}
\abs*{\frac{F(t)}{\norm{F(t)}}\paren*{\frac{g(t)}{\norm{f(t)}}}-1}=\abs*{\frac{\norm{g(t)}}{\norm{f(t)}}-1}<\beta_X\paren*{\frac{f(t)}{\norm{f(t)}},\frac{\eps}{4}}.
\end{equation*}
So by definition of $\beta_X$ we must have
\begin{equation*}
\frac{1}{2}\norm*{\frac{f(t)}{\norm{f(t)}}+\frac{g(t)}{\norm{f(t)}}}\leq 1-\beta_X\paren*{\frac{f(t)}{\norm{f(t)}},\frac{\eps}{4}}\leq1-\frac{1}{n_0}.
\end{equation*}
Once more by the definition of $B_1$ this implies
\begin{align*}
&\norm*{\frac{f(t)+g(t)}{2}}\leq\paren*{1-\frac{1}{n_0}}\norm{f(t)}\leq\frac{1-1/n_0}{2(1-\eta)}(\norm{f(t)}+\norm{g(t)}) \\
&=\frac{1}{2}(1-\alpha_1)(\norm{f(t)}+\norm{g(t)}),
\end{align*}
where $\alpha_1:=(1/n_0-\eta)(1-\eta)^{-1}>0$.\par
Now suppose that $t\in B_3\cap A_{n_0}\cap(S\sm N)$. Then
\begin{equation*}
\frac{F(t)}{\norm{F(t)}}\paren*{\frac{g(t)}{\norm{g(t)}}-\frac{f(t)}{\norm{g(t)}}}\geq\frac{\eps}{4},
\end{equation*}
consequently
\begin{align*}
&\frac{F(t)}{\norm{F(t)}}\paren*{\frac{g(t)}{\norm{g(t)}}-\frac{f(t)}{\norm{f(t)}}}\geq\frac{\eps}{4}+
\frac{F(t)}{\norm{F(t)}}\paren*{\frac{f(t)}{\norm{g(t)}}-\frac{f(t)}{\norm{f(t)}}} \\
&\geq\frac{\eps}{4}-\norm*{\frac{f(t)}{\norm{g(t)}}-\frac{f(t)}{\norm{f(t)}}}=\frac{\eps}{4}-\abs*{\frac{\norm{f(t)}}{\norm{g(t)}}-1}\geq\frac{\eps}{4}-\eta\geq\frac{\eps}{8}.
\end{align*}
Since $\norm{F(t)}\norm{g(t)}=F(t)(g(t))$ the definition of $\beta_X$ implies that
\begin{equation*}
\frac{1}{2}\norm*{\frac{f(t)}{\norm{f(t)}}+\frac{g(t)}{\norm{g(t)}}}\leq 1-\beta_X\paren*{\frac{f(t)}{\norm{f(t)}},\frac{\eps}{8}}\leq 1-\frac{1}{n_0},
\end{equation*}
where the latter inequality holds because of $t\in A_{n_0}$. It follows that
\begin{align*}
&\frac{1}{2}\norm*{\frac{f(t)}{\norm{f(t)}}+\frac{g(t)}{\norm{f(t)}}}\leq 1-\frac{1}{n_0}+\frac{1}{2}\norm*{\frac{g(t)}{\norm{f(t)}}-\frac{g(t)}{\norm{g(t)}}} \\
&=1-\frac{1}{n_0}+\frac{1}{2}\abs*{\frac{\norm{g(t)}}{\norm{f(t)}}-1}\leq 1-\frac{1}{n_0}+\frac{1}{2}\paren*{\frac{1}{1-\eta}-1}=1-\alpha_2,
\end{align*}
where $\alpha_2:=1/n_0-\eta(2-2\eta)^{-1}$ which by \eqref{eq:3.36} is greater than zero. Becuase of $\norm{f(t)}\leq\norm{g(t)}$ it follwos that
\begin{equation*}
\norm*{f(t)+g(t)}\leq (1-\alpha_2)(\norm{f(t)}+\norm{g(t)}).
\end{equation*}
So if we put $\alpha=\min\set*{\alpha_1,\alpha_2}$ and $P=B_2\cap A_{n_0}\cap(S\sm N)$, $Q=B_3\cap A_{n_0}\cap(S\sm N)$ then
\begin{equation}\label{eq:3.37}
\norm*{f(t)+g(t)}\leq (1-\alpha)(\norm{f(t)}+\norm{g(t)}) \ \ \forall t\in P\cup Q.
\end{equation}
Now we will show that if $i=2$ resp. $i=3$ then 
\begin{equation*}
\int_P \norm{F(t)}\norm{f(t)}\,\mathrm{d}\mu(t)\geq \frac{\eps}{64} \ \ \mathrm{reps.} \ \ \int_Q \norm{F(t)}\norm{f(t)}\,\mathrm{d}\mu(t)\geq \frac{\eps}{64}.
\end{equation*}
Let us first assume $i=2$, i.\,e.
\begin{equation*}
\int_{B_2} F(t)(g(t)-f(t))\,\mathrm{d}\mu(t)\geq \frac{\eps}{8}.
\end{equation*}
Since $\norm{f(t)}\geq\norm{g(t)}$ for $t\in B_2$ it follows that
\begin{equation*}
\int_{B_2} \norm{F(t)}\norm{f(t)}\,\mathrm{d}\mu(t)\geq\frac{\eps}{16}.
\end{equation*}
Because $N$ is a null set we have
\begin{align*}
&\int_P \norm{F(t)}\norm{f(t)}\,\mathrm{d}\mu(t)=\int_{B_2\cap A_{n_0}} \norm{F(t)}\norm{f(t)}\,\mathrm{d}\mu(t) \\
&=\int_{B_2} \norm{F(t)}\norm{f(t)}\,\mathrm{d}\mu(t)-\int_{B_2\sm A_{n_0}} \norm{F(t)}\norm{f(t)}\,\mathrm{d}\mu(t) \\
&\geq\frac{\eps}{16}-\int_{S\sm A_{n_0}} \norm{F(t)}\norm{f(t)}\,\mathrm{d}\mu(t)\geq\frac{\eps}{16}-\norm{\norm{f(\cdot)}\chi_{S\sm A_{n_0}}}_E \\
&\geq\frac{\eps}{16}-\frac{\eps}{64}\geq\frac{\eps}{64},
\end{align*}
where the second last inequality holds because of \eqref{eq:3.32}.\par
Now assume that $i=3$, i.\,e.
\begin{equation*}
\int_{B_3} F(t)(g(t)-f(t))\,\mathrm{d}\mu(t)\geq \frac{\eps}{8}.
\end{equation*}
It follows that
\begin{align*}
&\frac{\eps}{8}\leq \int_{B_3} \norm{F(t)}(\norm{g(t)}+\norm{f(t)})\,\mathrm{d}\mu(t) \\
&\leq \int_{B_3} \norm{F(t)}\paren*{1+\frac{1}{1-\eta}}\norm{f(t)}\,\mathrm{d}\mu(t)\leq 4\int_{B_3} \norm{F(t)}\norm{f(t)}\,\mathrm{d}\mu(t)
\end{align*}
and hence as before we get 
\begin{equation*}
\int_Q \norm{F(t)}\norm{f(t)}\,\mathrm{d}\mu(t)\geq \frac{\eps}{32}-\frac{\eps}{64}=\frac{\eps}{64}.
\end{equation*}
So if $i=2$ or $i=3$ then there is $R\in\set*{P,Q}$ such that
\begin{equation*}
\norm{\norm{f(\cdot)}\chi_R}_E\geq\int_R \norm{F(t)}\norm{f(t)}\,\mathrm{d}\mu(t)\geq \frac{\eps}{64}.
\end{equation*}
Put $h=\norm{f(\cdot)}(1-2\alpha\chi_R)$. Then $h\in B_E$ and moreover $\norm{\norm{f(\cdot)}-h}_E=
2\alpha\norm{\norm{f(\cdot)}\chi_R}_E\geq \alpha\eps/32$, hence 
\begin{equation*}
\norm{\norm{f(\cdot)}(1-\alpha\chi_R)}_E=\frac{1}{2}\norm{\norm{f(\cdot)}+h}_E\leq 1-\delta_E\paren*{\norm{f(\cdot)},\frac{\eps\alpha}{32}}.
\end{equation*}
We further have
\begin{align*}
&\norm*{\frac{f+g}{2}}_{E(X)}\leq \frac{1}{2}\norm{(\norm{f(\cdot)}+\norm{g(\cdot)})\chi_{S\sm R}+\norm{f(\cdot)+g(\cdot)}\chi_R}_E \\
&\stackrel{\eqref{eq:3.37}}{\leq} \frac{1}{2}\norm{(\norm{f(\cdot)}+\norm{g(\cdot)})\chi_{S\sm R}+(1-\alpha)(\norm{f(\cdot)}+\norm{g(\cdot)})\chi_R}_E \\
&\leq \frac{1}{2}\norm{\norm{g(\cdot)}+\norm{f(\cdot)}-\alpha\norm{f(\cdot)}\chi_R}_E\leq \frac{1}{2}+\frac{1}{2}\norm{\norm{f(\cdot)}(1-\alpha\chi_R)}_E \\
&\leq \frac{1}{2}+\frac{1}{2}\paren*{1-\delta_E\paren*{\norm{f(\cdot)},\frac{\eps\alpha}{32}}}=1-\frac{1}{2}\delta_E\paren*{\norm{f(\cdot)},\frac{\eps\alpha}{32}}.
\end{align*}
Altogether  we have shown that for
\begin{equation*}
\delta:=\min\set*{\frac{1}{2}\delta_E\paren*{\norm{f(\cdot)},\frac{\eps\alpha}{32}},\delta_E\paren*{\norm{f(\cdot)},\frac{\eps\eta}{16}}}>0
\end{equation*}
we have for every $g\in S_{E(X)}$ and every $l\in S_{E(X)^*}$ with $l(g)=1$ and $l(f)\leq 1-\eps$
\begin{equation*}
\norm*{\frac{f+g}{2}}_{E(X)}\leq 1-\delta.
\end{equation*}
By the aforementioneed characterisation of \ac{sluacs} spaces (\cite{hardtke}*{Proposition 2.1}) this implies that $E(X)$ is \ac{sluacs}.
\end{Proof}

Next we will have a look at the case of \ac{wuacs} spaces.
\begin{theorem}\label{thm:wuacs}
If $\mu$ is a $\sigma$-finite measure and $E$ is \ac{wuacs}, reflexive and has the Kadets-Klee property, then $E(X)$ 
is \ac{wuacs} whenever $X$ is \ac{wuacs}.
\end{theorem}

\begin{Proof}
Note that since $E$ is reflexive (or since it has the Kadets-Klee property), it is order continuous.\par
Let us take two sequences $(f_n)_{n\in \N}$ and $(g_n)_{n\in \N}$ in the unit sphere of $E(X)$ such 
that $\norm{f_n+g_n}_{E(X)}\to 2$ and a functional $l\in S_{E(X)^*}$, as usual represented by 
$[F]\in E^{\prime}(X^*,w^*)$, with $l(f_n)\to 1$.\par
As in the proof of Theorem \ref{thm:luacs} we find 
\begin{equation}\label{eq:3.38}
\lim_{n\to \infty}\int_S \norm{F(t)}\norm{f_n(t)}\,\mathrm{d}\mu(t)=1
\end{equation}
and by passing to a subsequence also
\begin{equation}\label{eq:3.39}
\lim_{n\to \infty}\paren*{\norm{F(t)}\norm{f_n(t)}-F(t)(f_n(t))}=0 \ \ \text{a.\,e.}
\end{equation}
It is also easy to see that
\begin{equation}\label{eq:3.40}
\lim_{n\to \infty}\norm*{\norm{f_n(\cdot)}+\norm{g_n(\cdot)}}_E=2
\end{equation}
and
\begin{equation}\label{eq:3.41}
\lim_{n\to \infty}\norm*{\norm{f_n(\cdot)}+\norm{g_n(\cdot)}+\norm{f_n(\cdot)+g_n(\cdot)}}_E=4.
\end{equation}
Since $E$ is \ac{wuacs} it follows from \eqref{eq:3.38} and \eqref{eq:3.40} that
\begin{equation}\label{eq:3.42}
\lim_{n\to \infty}\int_S \norm{F(t)}\norm{g_n(t)}\,\mathrm{d}\mu(t)=1.
\end{equation}
Again since $E$ is \ac{wuacs} and because of \eqref{eq:3.38}, \eqref{eq:3.40}, \eqref{eq:3.41} 
and \eqref{eq:3.42} we can deduce that
\begin{equation}\label{eq:3.43}
\lim_{n\to \infty}\int_S\norm{F(t)}\paren*{\norm{f_n(t)}+\norm{g_n(t)}-\norm{f_n(t)+g_n(t)}}\,\mathrm{d}\mu(t)=0.
\end{equation}
and hence we can pass to a further subsequence such that
\begin{equation}\label{eq:3.44}
\lim_{n\to \infty}\norm{F(t)}\paren*{\norm{f_n(t)}+\norm{g_n(t)}-\norm{f_n(t)+g_n(t)}}=0 \ \ \text{a.\,e.}
\end{equation}
By the reflexivity of $E$ we can pass once more to a subsequence such that $(\norm{f_n(\cdot)})_{n\in \N}$ and
$(\norm{g_n(\cdot)})_{n\in \N}$ are weakly convergent to $h_1\in B_E$ resp. $h_2\in B_E$. In view of \eqref{eq:3.38} 
and \eqref{eq:3.42} it follows that 
\begin{equation*}
\int_S\norm{F(t)}h_i(t)\,\mathrm{d}\mu(t)=1 \ \ \forall i\in\set*{1,2},
\end{equation*}
hence $\norm{h_1}_E=\norm{h_2}_E=1$ and moreover
\begin{equation}\label{eq:3.45}
\norm{h_1+h_2}_E=2.
\end{equation}
The fact that $E$ has the Kadets-Klee property implies that
\begin{equation*}
\norm*{\norm{f_n(\cdot)}-h_1}_E\to 0 \ \ \mathrm{and} \ \ \norm*{\norm{g_n(\cdot)}-h_2}_E\to 0
\end{equation*}
and thus by Lemma \ref{lemma:conv ae} we can, for the last time, pass to a subsequence such that
\begin{equation}\label{eq:3.46}
\lim_{n\to \infty}\norm{f_n(t)}=h_1(t) \ \ \text{and} \ \ \lim_{n\to \infty}\norm{g_n(t)}=h_2(t) \ \ \text{a.\,e.}
\end{equation}
Let $N_1$ resp. $N_2$ resp. $N_3$ denote the null sets on which the convergence statement from \eqref{eq:3.39} 
resp. \eqref{eq:3.44} resp. \eqref{eq:3.46} does not hold and put $N=N_1\cup N_2\cup N_3$ as well as
$B=\set*{t\in S\sm N:F(t)\neq 0 \ \mathrm{and} \ h_2(t)\neq 0}$ and $C=\set*{t\in B:h_1(t)=0}$.\par
Because of \eqref{eq:3.45} and since $E$ is in particular \ac{acs} we can show just as in the proof
of Proposition \ref{prop:acs} that $C$ is a null set.\par
The fact that $X$ is \ac{wuacs} together with Lemma \ref{aux lemma} easily implies that 
\begin{equation}\label{eq:3.47}
\lim_{n\to \infty}\paren*{\norm{F(t)}\norm{g_n(t)}-F(t)(g_n(t))}=0 \ \ \forall t\in S\sm (N\cup C).
\end{equation}
By the weak convergence of $(\norm{g_n(\cdot)})_{n\in \N}$ to $h_2$ we have
\begin{equation}\label{eq:3.48}
\lim_{n\to \infty}\int_A\norm{F(t)}\norm{g_n(t)}\,\mathrm{d}\mu(t)=\int_A\norm{F(t)}h_2(t)\,\mathrm{d}\mu(t) \ \ \forall A\in \A.
\end{equation}
Since $E$ is reflexive $E^{\prime}$ is order continuous and thus we can deduce as in the proof of Theorem \ref{thm:luacs+}, 
with the aid of Vitali's Lemma, \eqref{eq:3.48}, \eqref{eq:3.47} and the fact that $N\cup C$ is a null set, that
\begin{equation*}
\lim_{n\to \infty}\int_S\paren*{\norm{F(t)}\norm{g_n(t)}-F(t)(g_n(t))}\,\mathrm{d}\mu(t)=0.
\end{equation*}
Because of \eqref{eq:3.42} it follows that
\begin{equation*}
\lim_{n\to \infty}l(g_n)=\lim_{n\to \infty}\int_S F(t)(g_n(t))\,\mathrm{d}\mu(t)=1
\end{equation*}
and we are done.
\end{Proof}

If we combine the techniques of the proofs of Theorem \ref{thm:wuacs} and Theorem \ref{thm:luacs+} we can 
also obtain another result concerning $\mathrm{luacs}^+$ spaces (we omit the details).
\begin{theorem}\label{thm:luacs+2}
If $\mu$ is a $\sigma$-finite measure and $E$ is $\text{luacs}^+$, reflexive and has the Kadets-Klee property, 
then $E(X)$ is $\text{luacs}^+$ whenever $X$ is $\text{luacs}^+$.
\end{theorem}

It is further possible to obtain another sufficient condition for $E(X)$ to be \ac{sluacs}.
\begin{theorem}\label{thm:sluacs2}
If $\mu$ is a $\sigma$-finite measure and $E$ is $\text{sluacs}^+$ and reflexive and both $E$ and $E^*$ 
have the Kadets-Klee property, then $E(X)$ is \ac{sluacs} whenever $X$ is \ac{sluacs}.
\end{theorem}

\begin{Proof}
Let $(f_n)_{n\in \N}$ be a sequence in $S_{E(X)}$ and $f\in S_{E(X)}$ such that we have $\norm{f_n+f}_{E(X)}\to 2$.
Also, let $(l_n)_{n\in \N}$ be a sequence in $S_{E(X)^*}$ such that $l_n(f_n)\to 1$. If we represent each 
$l_n$ by $[F_n]\in E^{\prime}(X^*,w^*)$ we can obtain as usual
\begin{equation}\label{eq:3.49}
\lim_{n\to \infty}\int_S \norm{F_n(t)}\norm{f_n(t)}\,\mathrm{d}\mu(t)=1
\end{equation}
and by passing to a subsequence also 
\begin{equation}\label{eq:3.50}
\lim_{n\to \infty}\paren*{\norm{F_n(t)}\norm{f_n(t)}-F_n(t)(f_n(t))}=0 \ \ \text{a.\,e.}
\end{equation}
as well as
\begin{align}
&\lim_{n\to \infty}\norm*{\norm{f_n(\cdot)}+\norm{f(\cdot)}}_E=2,\label{eq:3.51} \\
&\lim_{n\to \infty}\norm*{\norm{f_n(\cdot)+f(\cdot)}+\norm{f_n(\cdot)}+\norm{f(\cdot)}}_E=4,\label{eq:3.52} \\
&\lim_{n\to \infty}\norm*{\norm{f_n(\cdot)+f(\cdot)}+\norm{f_n(\cdot)}}_E=3,\label{eq.3.53} \\
&\lim_{n\to \infty}\norm*{\norm{f_n(\cdot)+f(\cdot)}+\norm{f(\cdot)}}_E=3,\label{eq:3.54} \\
&\lim_{n\to \infty}\norm*{\norm{f_n(\cdot)+f(\cdot)}+\norm{f_n(\cdot)}+3\norm{f(\cdot)}}_E=6.\label{eq:3.55}
\end{align}
Using the fact that $E$ is $\mathrm{sluacs}^+$ we can conclude that
\begin{equation}\label{eq:3.56}
\lim_{n\to \infty}\int_S \norm{F_n(t)}\norm{f(t)}\,\mathrm{d}\mu(t)=1
\end{equation}
and
\begin{equation}\label{eq:3.57}
\lim_{n\to \infty}\int_S \norm{F_n(t)}\paren*{\norm{f_n(t)}+\norm{f(t)}-\norm{f_n(t)+f(t)}}\,\mathrm{d}\mu(t)=0.
\end{equation}
So we can pass to another subsequence such that
\begin{equation}\label{eq:3.58}
\lim_{n\to \infty}\norm{F_n(t)}\paren*{\norm{f_n(t)}+\norm{f(t)}-\norm{f_n(t)+f(t)}}=0 \ \ \text{a.\,e.}
\end{equation}
Since $E$ (and hence also $E^*$) is reflexive we may assume without loss of generality that $(\norm{f_n(\cdot)})_{n\in \N}$
is weakly convergent to some $h\in B_E$ and that $(\norm{F_n(\cdot)})_{n\in \N}$ is weakly convergent to some $g\in B_{E^*}=B_{E^{\prime}}$.\par
It follows from \eqref{eq:3.56} that 
\begin{equation}\label{eq:3.59}
\int_S g(t)\norm{f(t)}\,\mathrm{d}\mu(t)=1
\end{equation}
and hence $g\in S_{E^*}$. Because of \eqref{eq:3.59}, \eqref{eq:3.51} and the fact that $E$ is $\mathrm{sluacs}^+$
we get that
\begin{equation*}
\lim_{n\to \infty}\int_S g(t)\norm{f_n(t)}\,\mathrm{d}\mu(t)=1
\end{equation*}
and consequently
\begin{equation}\label{eq:3.60}
\int_S g(t)h(t)\,\mathrm{d}\mu(t)=1,
\end{equation}
whence $h\in S_E$. Since both $E$ and $E^*$ have the Kadets-Klee property it follows that
\begin{equation}\label{eq:3.61}
\norm*{\norm{f_n(\cdot)}-h}_E\to 0 \ \ \mathrm{and} \ \ \norm*{\norm{F_n(\cdot)}-g}_{E^{\prime}}\to 0.
\end{equation}
Thus we can pass once more to a subsequence such that
\begin{equation}\label{eq:3.62}
\lim_{n\to \infty}\norm{f_n(t)}=h(t) \ \ \mathrm{and} \ \ \lim_{n\to \infty}\norm{F_n(t)}=g(t) \ \ \text{a.\,e.}
\end{equation}
Combining \eqref{eq:3.60} and \eqref{eq:3.59} we also obtain
\begin{equation}\label{eq:3.63}
\norm{h+\norm{f(\cdot)}}_E=2.
\end{equation}
Let $N$ be a null set such that the convergence statements of \eqref{eq:3.50}, \eqref{eq:3.58} and \eqref{eq:3.62}
hold for every $t\in S\sm N$ and put $B=\set*{t\in S\sm N:g(t)\neq 0 \ \mathrm{and} f(t)\neq 0}$ as well as 
$C=\set*{t\in B:h(t)=0}$.\par
Similar to the arguments in the proof of Theorem \ref{thm:wuacs} one can see that $C$ is a null set
and then, using the fact that $X$ is \ac{sluacs}, deduce that
\begin{equation*}
\lim_{n\to \infty}\paren*{\norm{F_n(t)}\norm{f(t)}-F_n(t)(f(t))}=0 \ \ \text{a.\,e.}
\end{equation*}
By our usual method based on Vitali's Lemma we can conclude that for every $A\in \A$ with $\mu(A)<\infty$ we have
\begin{equation}\label{eq:3.64}
\lim_{n\to \infty}\int_A\paren*{\norm{F_n(t)}\norm{f(t)}-F_n(t)(f(t))}\,\mathrm{d}\mu(t)=0.
\end{equation}
Now we fix an increasing sequence $(A_m)_{m\in \N}$ in $\A$ such that $\mu(A_m)<\infty$ for all $m\in \N$ and $\bigcup_{m=1}^{\infty}A_m=S$.
The order continuity of $E$ implies $\norm*{\norm{f(\cdot)}\chi_{S\sm A_m}}_E\to 0$. Analogous to the argument at the end of the proof of
Theorem \ref{thm:luacs+} this together with \eqref{eq:3.64} leads to
\begin{equation*}
\lim_{n\to \infty}\int_S\paren*{\norm{F_n(t)}\norm{f(t)}-F_n(t)(f(t))}\,\mathrm{d}\mu(t)=0.
\end{equation*}
Taking into account \eqref{eq:3.56} we arrive at
\begin{equation*}
\lim_{n\to \infty}l_n(f)=\lim_{n\to \infty}\int_S F_n(t)(f(t))\,\mathrm{d}\mu(t)=1
\end{equation*}
and the proof is finished.
\end{Proof}

Next we will consider sufficient conditions for a K\"othe-Bochner space to be $\mathrm{sluacs}^+$
(recall that a dual Banach space $X^*$ is said to have the Kadets-Klee* property if it fulfils the
definition of the Kadets-Klee property with weak- replaced by weak*-convergence).
\begin{theorem}\label{thm:sluacs+}
Let $E$ be a K\"othe function space over the complete $\sigma$-finite measure space $(S,\A,\mu)$
and let $X$ be an $\text{sluacs}^+$ Banach space. If $E^*$ has the Kadets-Klee* property and in 
addition
\begin{enumerate}[\upshape(a)]
\item $E$ is $\text{sluacs}^+$, reflexive and has the Kadets-Klee property or
\item $E$ is \ac{LUR} and $B_{E^*}$ is weak*-sequentially compact,
\end{enumerate}
then $E(X)$ is $\text{sluacs}^+$.
\end{theorem}

\begin{Proof}
By the Theorems \ref{thm:sluacs} and \ref{thm:sluacs2} we already know that $E(X)$ is in both cases \ac{sluacs}.
Note also that in both cases $E$ is order continuous. Now take a sequence $(f_n)_{n\in \N}$ in $S_{E(X)}$ and $f\in S_{E(X)}$ 
such that $\norm{f_n+f}_{E(X)}\to 2$ and let $(l_n)_{n\in \N}$ be a sequence in $S_{E(X)^*}$ such that $l_n(f)\to 1$.
If we represent each $l_n$ by $[F_n]\in E^{\prime}(X^*,w^*)$ we can obtain as usual
\begin{equation}\label{eq:3.65}
\lim_{n\to \infty}\int_S \norm{F_n(t)}\norm{f(t)}\,\mathrm{d}\mu(t)=1
\end{equation}
and by passing to a subsequence also 
\begin{equation}\label{eq:3.66}
\lim_{n\to \infty}\paren*{\norm{F_n(t)}\norm{f(t)}-F_n(t)(f(t))}=0 \ \ \text{a.\,e.}
\end{equation}
as well as
\begin{align}
&\lim_{n\to \infty}\norm*{\norm{f_n(\cdot)}+\norm{f(\cdot)}}_E=2,\label{eq:3.67} \\
&\lim_{n\to \infty}\norm*{\norm{f_n(\cdot)+f(\cdot)}+\norm{f_n(\cdot)}+\norm{f(\cdot)}}_E=4,\label{eq:3.68} \\
&\lim_{n\to \infty}\norm*{\norm{f_n(\cdot)+f(\cdot)}+\norm{f_n(\cdot)}}_E=3,\label{eq.3.69} \\
&\lim_{n\to \infty}\norm*{\norm{f_n(\cdot)+f(\cdot)}+\norm{f(\cdot)}}_E=3,\label{eq:3.70} \\
&\lim_{n\to \infty}\norm*{\norm{f_n(\cdot)+f(\cdot)}+\norm{f_n(\cdot)}+3\norm{f(\cdot)}}_E=6.\label{eq:3.71}
\end{align}
Since $E$ is $\mathrm{sluacs}^+$ it follows that
\begin{equation}\label{eq:3.72}
\lim_{n\to \infty}\int_S \norm{F_n(t)}\norm{f_n(t)}\,\mathrm{d}\mu(t)=1
\end{equation}
and
\begin{equation}\label{eq:3.73}
\lim_{n\to \infty}\int_S \norm{F_n(t)}\paren*{\norm{f_n(t)}+\norm{f(t)}-\norm{f_n(t)+f(t)}}\,\mathrm{d}\mu(t)=0,
\end{equation}
so that by passing to another subsequence we can assume
\begin{equation}\label{eq:3.74}
\lim_{n\to \infty}\norm{F_n(t)}\paren*{\norm{f_n(t)}+\norm{f(t)}-\norm{f_n(t)+f(t)}}=0 \ \ \text{a.\,e.}
\end{equation}
In both cases (a) and (b) the dual unit ball $B_{E^*}$ is weak*-sequentially compact so that we can also assume
the weak*-convergence of $(\norm{F_n(\cdot)})_{n\in \N}$ to some $g\in B_{E^*}$. It follows from \eqref{eq:3.65} that
\begin{equation}\label{eq:3.75}
\int_S g(t)\norm{f(t)}\,\mathrm{d}\mu(t)=1
\end{equation}
and hence $\norm{g}_{E^{\prime}}=1$. Since $E^*$ has the Kadets-Klee* property we get that
\begin{equation}\label{eq:3.76}
\norm*{\norm{F_n(\cdot)}-g}_{E^{\prime}}\to 0
\end{equation}
and thus we can, by passing to yet another subsequence, assume that
\begin{equation}\label{eq:3.77}
\lim_{n\to \infty}\norm{F_n(t)}=g(t) \ \ \text{a.\,e.}
\end{equation}
Next we claim that there is an $h\in S_E$ such that
\begin{equation}\label{eq:3.78}
\int_S g(t)h(t)\,\mathrm{d}\mu(t)=1
\end{equation}
and, after passing to a subsequence once more,
\begin{equation}\label{eq:3.79}
\norm*{\norm{f_n(\cdot)}-h}_E\to 0.
\end{equation}
For in the case (b) $E$ is \ac{LUR} and thus by \eqref{eq:3.67} and \eqref{eq:3.75} we can take $h=\norm{f(\cdot)}$.
In the case (a) $E$ is reflexive and hence we can assume that $(\norm{f_n(\cdot)})_{n\in \N}$ is weakly convergent
to some $h\in B_E$. Then \eqref{eq:3.78} follows from \eqref{eq:3.76} and \eqref{eq:3.72}. This also implies $\norm{h}_E=1$
and by the Kadets-Klee property of $E$ we have \eqref{eq:3.79}.\par
By \eqref{eq:3.79} we may assume that
\begin{equation}\label{eq:3.80}
\lim_{hn\to \infty}\norm{f_n(t)}=h(t) \ \ \text{a.\,e.}
\end{equation}
Note that \eqref{eq:3.75} and \eqref{eq:3.78} imply that $\norm*{\norm{f(\cdot)}+h}_E=2$. Using all this and the fact that
$X$ is $\mathrm{sluacs}^+$ one can first prove, analogously to the arguments in the proof of Theorem \ref{thm:sluacs2}, that
\begin{equation}\label{eq:3.81}
\lim_{n\to \infty}\paren*{\norm{F_n(t)}\norm{f_n(t)}-F_n(t)(f_n(t))}=0 \ \ \text{a.\,e.}
\end{equation}
and then
\begin{equation}\label{eq:3.82}
\lim_{n\to \infty}\int_A\paren*{\norm{F_n(t)}\norm{f_n(t)}-F_n(t)(f_n(t))}\,\mathrm{d}\mu(t)=0
\end{equation}
for every $A\in \A$ with $\mu(A)<\infty$.\par
Let us now fix a sequence $(A_m)_{m\in \N}$ in $\A$ as in the proof of Theorem \ref{thm:sluacs2}. The order continuity of $E$ implies
$\norm{\norm{f(\cdot)}\chi_{S\sm A_m}}_E\to 0$.\par
Let $\eps>0$ be arbitrary. Since $E$ is $\mathrm{sluacs}^+$ there exists a $\delta>0$ such that for all $b\in S_E$
and all $l\in B_{E^*}$ with $\norm*{b+\norm{f(\cdot)}}_E\geq 2(1-\delta)$ and $l(\norm{f(\cdot)})\geq 1-\delta$ 
one has $l(b)\geq 1-\eps$.\par
Fix $m_0\in \N$ with $\norm{\norm{f(\cdot)}\chi_{S\sm A_{m_0}}}_E\leq \delta/2$. Because of \eqref{eq:3.67}, \eqref{eq:3.65} and \eqref{eq:3.82}
there is an $N\in \N$ such that for all $n\geq N$ the inequalities
\begin{align*}
&\norm*{\norm{f_n(\cdot)}+\norm{f(\cdot)}}_E\geq 2(1-\delta), \\
&\int_S\norm{F_n(t)}\norm{f(t)}\,\mathrm{d}\mu(t)\geq 1-\frac{\delta}{2}, \\
&\int_{A_{m_0}}\paren*{\norm{F_n(t)}\norm{f_n(t)}-F_n(t)(f_n(t))}\,\mathrm{d}\mu(t)\leq\eps
\end{align*}
hold.\par
It follows that for every $n\geq N$ we have
\begin{align*}
&\abs*{\int_{A_{m_0}}\norm{F_n(t)}\norm{f(t)}\,\mathrm{d}\mu(t)-1} \\
&\leq\norm{\norm{f(\cdot)}\chi_{S\sm A_{m_0}}}_E+\abs*{\int_S\norm{F_n(t)}\norm{f(t)}\,\mathrm{d}\mu(t)-1}\leq\delta
\end{align*}
and hence by the choice of $\delta$ 
\begin{equation*}
\int_{A_{m_0}}\norm{F_n(t)}\norm{f_n(t)}\,\mathrm{d}\mu(t)\geq 1-\eps.
\end{equation*}
Consequently, for every $n\geq N$ we have
\begin{align*}
&\int_S\paren*{\norm{F_n(t)}\norm{f_n(t)}-F_n(t)(f_n(t))}\,\mathrm{d}\mu(t)\leq \\
&\leq\eps+\int_{S\sm A_{m_0}}\paren*{\norm{F_n(t)}\norm{f_n(t)}-F_n(t)(f_n(t))}\,\mathrm{d}\mu(t) \\
&\leq\eps+2\int_{S\sm A_{m_0}}\norm{F_n(t)}\norm{f_n(t)}\,\mathrm{d}\mu(t)\leq \eps+2(1-(1-\eps))=3\eps.
\end{align*}
Thus we have shown 
\begin{equation*}
\lim_{n\to \infty}\int_S\paren*{\norm{F_n(t)}\norm{f_n(t)}-F_n(t)(f_n(t))}\,\mathrm{d}\mu(t)=0.
\end{equation*}
Together with \eqref{eq:3.72} it follows $l_n(f_n)\to 1$, as desired.
\end{Proof}

Now let us treat the case of \ac{uacs} spaces. In analogy to \cite{hardtke}*{Definition 3.15}
we say that an order continuous K\"othe function space $E$ has property $(u^+)$ if for every $\eps>0$
there is some $\delta>0$ such that for all $f,g\in S_E$ and every $h\in S_{E^{\prime}}$ we have
\begin{equation*}
\norm{f+g}_E\geq2(1-\delta) \ \ \mathrm{and} \ \ \int_S fh\,\mathrm{d}\mu=1 \ \ \Rightarrow \ \ \int_S \abs{h}\abs{f-g}\,\mathrm{d}\mu\leq\eps.
\end{equation*}
This property certainly implies that $E$ is \ac{uacs}. Every \ac{UR} space has property $(u^+)$.
The following theorem holds. Its proof is completely analogous to the one of \cite{hardtke}*{Theorem 3.16}
(which is a modification of the proof of \cite{day}*{Theorem 3}) but we will explicitly give it here, for 
the readers convenience.
\begin{theorem}\label{thm:uacs u+}
If $E$ is an order continuous K\"othe function space with the property $(u^+)$ (in particular, if $E$ is \ac{UR})
and $X$ is a \ac{uacs} Banach space then $E(X)$ is also \ac{uacs}.
\end{theorem}

\begin{Proof}
Let $0<\eps\leq 2$ be arbitrary. Since $E$ is in particular \ac{uacs} there is a number $\eta>0$ such that
for all functions $a,b\in B_E$ and every functional $l\in B_{E^*}$ with $l(a)=1$ one has
\begin{equation}\label{eq:3.83}
l(b)<1-\frac{\eps}{4}\delta_{\mathrm{uacs}}^X(\eps/2) \ \Rightarrow \ \norm{a+b}_E\leq 2(1-\eta).
\end{equation}
Now let $f,g\in S_{E(X)}$ such that $\norm{f(t)}=\norm{g(t)}$ a.\,e. and let $L\in E(X)^*$ such that
$L(f)=1$ and $L(g)<1-\eps$. We claim that $\norm{f+g}_{E(X)}\leq 2(1-\eta)$.\par
Let $L$ be represented by $[F]\in E^{\prime}(X^*,w^*)$ and put $\beta=\norm{g(\cdot)}$, $\nu=\norm{F(\cdot)}$.
Define $\gamma$ by $\gamma(t)=\nu(t)\beta(t)-F(t)(g(t))$. Note that $\gamma$ is measurable and
\begin{equation}\label{eq:3.84}
0\leq\gamma(t)\leq2\nu(t)\beta(t) \ \ \forall t\in S.
\end{equation}
As before we can deduce from $L(f)=1$ that
\begin{equation}\label{eq:3.85}
\int_S\norm{F(t)}\norm{f(t)}\,\mathrm{d}\mu(t)=1
\end{equation}
and $F(t)(f(t))=\norm{F(t)}\norm{f(t)}$ a.\,e., hence
\begin{equation}\label{eq:3.86}
F(t)(f(t))=\nu(t)\beta(t) \ \ \text{a.\,e.}
\end{equation}
Next we define
\begin{equation*}
\alpha(t)=\begin{cases}
&\frac{1}{2}\delta_{\mathrm{uacs}}^X\paren*{\frac{\gamma(t)}{\nu(t)\beta(t)}} \ \text{if} \ 0<\gamma(t)<\nu(t)\beta(t) \\
&0 \ \text{if} \ \gamma(t)=0 \\
&\frac{1}{2}\delta_{\mathrm{uacs}}^X(1) \ \text{otherwise.}
\end{cases}
\end{equation*}
Note that since $\delta_{\mathrm{uacs}}^X$ is continuous on $(0,1)$ (see \cite{dhompongsa0}*{Lemma 3.10} or \cite{hardtke}*{Lemma 2.11}), 
the function $\alpha$ is measurable. Using \eqref{eq:3.86} it is easy to see that
\begin{equation}\label{eq:3.87}
\norm{f(t)+g(t)}\leq 2(1-\alpha(t))\beta(t) \ \ \text{a.\,e.}
\end{equation}
By \eqref{eq:3.84} and \eqref{eq:3.85} we have $\int_S\gamma(t)\,\mathrm{d}\mu(t)\leq 2$. Furthermore, we also have
\begin{equation*}
\eps<1-L(g)=L(f-g)=\int_S F(t)(f(t)-g(t))\,\mathrm{d}\mu(t)\leq\int_S \gamma(t)\,\mathrm{d}\mu(t),
\end{equation*}
thus
\begin{equation}\label{eq:3.88}
\eps<\int_S \gamma(t)\,\mathrm{d}\mu(t)\leq 2.
\end{equation}
Now put $A=\set*{t\in S:2\gamma(t)>\eps\nu(t)\beta(t)}$ and $B=S\sm A$. We then have (because of \eqref{eq:3.85})
\begin{equation*}
\int_B \gamma(t)\,\mathrm{d}\mu(t)\leq\frac{\eps}{2}\int_B \nu(t)\beta(t)\,\mathrm{d}\mu(t)\leq\frac{\eps}{2}\int_S \nu(t)\beta(t)\,\mathrm{d}\mu(t)=\frac{\eps}{2}.
\end{equation*}
Together with \eqref{eq:3.88} it follows that 
\begin{equation*}
\int_A \gamma(t)\,\mathrm{d}\mu(t)>\eps-\frac{\eps}{2}=\frac{\eps}{2}.
\end{equation*}
Taking into account \eqref{eq:3.84} we get
\begin{equation}\label{eq:3.89}
\int_A \nu(t)\beta(t)\,\mathrm{d}\mu(t)>\frac{\eps}{4}.
\end{equation}
Next we define $h=\beta\chi_B$ and $h^{\prime}=\beta\chi_A$ as well as $h^{\prime\prime}=(1-\delta_{\mathrm{uacs}}^X(\eps/2))h^{\prime}$.
Then $\norm{h+h^{\prime\prime}}_E\leq\norm{h+h^{\prime}}_E=\norm{\beta}_E=1$. Let $l$ be the functional on $E$ represented by $\nu=\norm{F(\cdot)}$.
We have $l(h+h^{\prime})=l(\beta)=1$ (by \eqref{eq:3.85}) and further, by \eqref{eq:3.89},
\begin{equation*}
l(h+h^{\prime\prime})=1-\delta_{\mathrm{uacs}}^X(\eps/2)l(h^{\prime})=1-\int_A \nu(t)\beta(t)\,\mathrm{d}\mu(t)<1-\frac{\eps}{4}\delta_{\mathrm{uacs}}^X(\eps/2).
\end{equation*}
So by our choice of $\eta$ we get $\norm{2h+h^{\prime}+h^{\prime\prime}}_E\leq 2(1-\eta)$, i.\,e.
\begin{equation}\label{eq:3.90}
\norm*{h+\paren*{1-\frac{1}{2}\delta_{\mathrm{uacs}}^X(\eps/2)}h^{\prime}}_E\leq1-\eta.
\end{equation}
By monotonicity of $\delta_{\mathrm{uacs}}^X$ we have
\begin{equation}\label{eq:3.91}
\alpha(t)\geq\frac{1}{2}\delta_{\mathrm{uacs}}^X(\eps/2) \ \ \forall t\in A.
\end{equation}
Using \eqref{eq:3.87}, \eqref{eq:3.91} and \eqref{eq:3.90} we obtain
\begin{align*}
&\norm{f+g}_{E(X)}=\norm{\norm{f(\cdot)+g(\cdot)}}_E\leq 2\norm{(1-\alpha)\beta}_E \\
&\leq 2\norm{(1-2^{-1}\delta_{\mathrm{uacs}}^X(\eps/2))h^{\prime}+h}_E\leq 2(1-\eta).
\end{align*}
The first step of the proof is completed. Next we wish to remove the restriction $\norm{f(\cdot)}=\norm{g(\cdot)}$ a.\,e.
So let again $0<\eps\leq 2$ be arbitrary and choose $\eta$ as above but corresponding to the value $\eps/2$. Take $0<\omega<2\eta/3$.\par
Since $E$ is \ac{uacs} we may find $\tau>0$ such that for all $a,b\in B_E$ and every $l\in B_{E^*}$ we have
\begin{equation}\label{eq:3.92}
l(a)\geq 1-\tau \ \mathrm{and} \ \norm{a+b}_E\geq 2(1-\tau) \ \Rightarrow \ l(b)\geq 1-\omega.
\end{equation}
Next we fix $0<\rho<\min\set*{\eps/2,2\tau,\omega}$ and find a number $\tilde{\tau}$ to the value $\rho$ according to the definition 
of the property $(u^+)$ of $E$. Finally, let $0<\xi<\min\set*{\tau,\tilde{\tau}}$.\par
Let $f,g\in S_{E(X)}$ be arbitrary and $L\in S_{E(X)^*}$ (as usually represented by $F$) such that $L(f)=1$ and $\norm{f+g}_{E(X)}\geq 2(1-\xi)$. 
We are going to prove that $L(g)>1-\eps$, thus showing that $E(X)$ is \ac{uacs}.\par
To this end, we define $z:S \rightarrow X$ by
\begin{equation*}
z(t)=\begin{cases}
&\frac{\norm{f(t)}}{\norm{g(t)}}g(t) \ \text{if} \ g(t)\neq 0 \\
&f(t) \ \text{if} \ g(t)=0.
\end{cases}
\end{equation*}
Then $z$ is Bochner-measurable and $\norm{z(t)}=\norm{f(t)}$ for all $t\in S$ (hence $z\in E(X)$). Furthermore,
\begin{equation}\label{eq:3.93}
\norm{z(t)-g(t)}=\abs{\norm{f(t)}-\norm{g(t)}} \ \ \forall t\in S.
\end{equation}
As before we have
\begin{equation}\label{eq:3.94}
\int_S \norm{F(t)}\norm{f(t)}\,\mathrm{d}\mu(t)=1.
\end{equation}
Also,
\begin{equation*}
2(1-\tilde{\tau})\leq 2(1-\xi)\leq\norm{f+g}_{E(X)}\leq\norm{\norm{f(\cdot)}+\norm{g(\cdot)}}_E,
\end{equation*}
so the choice of $\tilde{\tau}$ together with \eqref{eq:3.93} implies
\begin{equation}\label{eq:3.95}
\int_S \norm{F(t)}\norm{z(t)-g(t)}\,\mathrm{d}\mu(t)\leq\rho.
\end{equation}
Next we observe that
\begin{equation*}
\norm{\norm{f(\cdot)}+\norm{g(\cdot)}+\norm{f(\cdot)+g(\cdot)}}_E\geq2\norm{f+g}_{E(X)}\geq 4(1-\xi)\geq 4(1-\tau)
\end{equation*}
and (because of \eqref{eq:3.94} and \eqref{eq:3.95})
\begin{align*}
&\int_S\norm{F(t)}(\norm{f(t)}+\norm{g(t)})\,\mathrm{d}\mu(t)=1+\int_S \norm{F(t)}\norm{g(t)}\,\mathrm{d}\mu(t) \\
&\geq 1+\int_S \norm{F(t)}\norm{f(t)}\,\mathrm{d}\mu(t)-\int_S \norm{F(t)}\abs{\norm{f(t)}-\norm{g(t)}}\,\mathrm{d}\mu(t) \\
&=2-\int_S \norm{F(t)}\abs{\norm{f(t)}-\norm{g(t)}}\,\mathrm{d}\mu(t)\geq 2-\rho\geq 2(1-\tau).
\end{align*}
So \eqref{eq:3.92} implies 
\begin{equation}\label{eq:3.96}
\int_S \norm{F(t)}\norm{f(t)+g(t)}\,\mathrm{d}\mu(t)\geq 2(1-\omega).
\end{equation}
Using \eqref{eq:3.95} and \eqref{eq:3.96} we can conclude
\begin{align*}
&\norm{f+z}_{E(X)}\geq\int_S \norm{F(t)}\norm{f(t)+z(t)}\,\mathrm{d}\mu(t) \\
&\geq\int_S \norm{F(t)}\norm{f(t)+g(t)}\,\mathrm{d}\mu(t)-\int_S \norm{F(t)}\norm{g(t)-z(t)}\,\mathrm{d}\mu(t) \\
&\geq 2(1-\omega)-\rho>2(1-\eta).
\end{align*}
By the choice of $\eta$ this implies $L(z)\geq 1-\eps/2$. But by \eqref{eq:3.95} we also have $\abs{L(g)-L(z)}\leq\rho$, hence
$L(g)\geq L(z)-\rho\geq 1-\eps/2-\rho>1-\eps$.
\end{Proof}

The above theorem admits the following corollary.
\begin{corollary}\label{cor:US}
If $E$ is a \ac{US} K\"othe function space and $X$ is a \ac{uacs} Banach space then $E(X)$ is also \ac{uacs}.
\end{corollary}

\begin{Proof}
Since \ac{uacs} is a self-dual property (cf. \cite{hardtke}*{Corollary 2.13}) $X^*$ is also \ac{uacs} and 
since $E$ is \ac{US} we have that $E^*=E^{\prime}$ is \ac{UR}  (cf. \cite{fabian}*{Theorem 9.10}). So by the 
previous theorem $E^{\prime}(X^*)$ is \ac{uacs}. But as a \ac{uacs} space $X^*$ is reflexive and hence it has 
the Radon-Nikod\'ym property. It follows from the general theory in \cite{bukhvalov} that in this case $E(X)^*$ 
is isometrically isomorphic to $E^{\prime}(X^*)$, so $E(X)^*$ and hence also $E(X)$ is \ac{uacs}.
\end{Proof}

Finally, we consider some midpoint version of \ac{luacs} and \ac{sluacs} spaces. Let us first recall the
following well-known notions: a Banach space $X$ is said to be {\em \ac{MLUR}} if for any two sequences $(x_n)_{n\in \N}$
and $(y_n)_{n\in \N}$ in $S_X$ and every $x\in S_X$ we have
\begin{equation*}
\norm*{x-\frac{x_n+y_n}{2}}\to 0 \ \Rightarrow \ \norm{x_n-y_n}\to 0.
\end{equation*}
$X$ is called {\em \ac{WMLUR}} if it satisfies the above condition with $\norm{x_n-y_n}\to 0$ replaced 
by $x_n-y_n\xrightarrow{\sigma} 0$, where the symbol $\xrightarrow{\sigma}$ denotes the convergence in 
the weak topology of $X$. The notion of \ac{MLUR} spaces was originally introduced by Anderson in \cite{anderson}.\par
In \cite{hardtke} the author introduced the following analogous definitions.
\begin{definition}
Let $X$ be a Banach space.
  \begin{enumerate}[(i)]
  \item The space $X$ is said to be {\em \ac{mluacs}} if for any two sequences $(x_n)_{n\in \N}$ 
  and $(y_n)_{n\in \N}$ in $S_X$, every $x\in S_X$ and every $x^*\in S_{X^*}$ we have that
  \begin{equation*}
  \norm*{x-\frac{x_n+y_n}{2}}\to 0 \ \mathrm{and} \ x^*(x_n)\to 1 \ \Rightarrow \ x^*(y_n)\to 1.
  \end{equation*}
  \item The space $X$ is called {\em \ac{msluacs}} if for any two sequences $(x_n)_{n\in \N}$ and 
  $(y_n)_{n\in \N}$ in $S_X$, every $x\in S_X$ and every sequence $(x_n^*)_{n\in \N}$ in $S_{X^*}$ we 
  have that
  \begin{equation*}
  \norm*{x-\frac{x_n+y_n}{2}}\to 0 \ \mathrm{and} \ x_n^*(x_n)\to 1 \ \Rightarrow \ x_n^*(y_n)\to 1.
  \end{equation*}
  \end{enumerate}
\end{definition}

The chart below summarises the obvious implications. No other implications are true in general (see the examples in \cite{hardtke}).
\begin{figure}[H]
\begin{center}
  \begin{tikzpicture}
  \node (LUR) at (-2,0) {LUR};
  \node (MLUR) at (0,1) {MLUR};
  \node (WLUR) at (0,-1) {WLUR};
  \node (WMLUR) at (2,0) {WMLUR};
  \node (R) at (4,0) {R};
  \node (sluacs) at (-2,-1) {sluacs};
  \node (msluacs) at (0,0) {msluacs};
  \node (luacs) at (0,-2) {luacs};
  \node (mluacs) at (2,-1) {mluacs};
  \node (acs) at (4,-1) {acs};
  \draw[->] (LUR)--(MLUR);
  \draw[->] (LUR)--(WLUR);
  \draw[->] (WLUR)--(WMLUR);
  \draw[->] (MLUR)--(WMLUR);
  \draw[->] (WMLUR)--(R);
  \draw[->] (sluacs)--(msluacs);
  \draw[->] (sluacs)--(luacs);
  \draw[->] (luacs)--(mluacs);
  \draw[->] (msluacs)--(mluacs);
  \draw[->] (mluacs)--(acs);
  \draw[->] (LUR)--(sluacs);
  \draw[->] (MLUR)--(msluacs);
  \draw[->] (WLUR)--(luacs);
  \draw[->] (WMLUR)--(mluacs);
  \draw[->] (R)--(acs);
  \end{tikzpicture}
\end{center}
\CAP\label{fig:4}
\end{figure}

Concerning the properties \ac{msluacs} and \ac{mluacs} for K\"othe-Bochner spaces we have the following result.
\begin{theorem}\label{thm:midpoint}
Let $E$ be an \ac{MLUR} K\"othe function space over a complete $\sigma$-finite measure space and $X$ a Banach space.
If $X$ is \ac{mluacs}, then so is $E(X)$. If $X$ is \ac{msluacs} and in addition $E^*$ has the Kadets-Klee* property
and $B_{E^*}$ is weak*-sequentially compact, then $E(X)$ is also \ac{msluacs}.
\end{theorem}

\begin{Proof}
Let us first recall that $\ell^{\infty}$ has no equivalent \ac{MLUR} norm (cf. \cite{lin}*{Theorem 2.1.5}) and so
by \cite{lin}*{Propositions 3.1.4 and 3.1.5} (and since every K\"othe function space is $\sigma$-order complete) $E$ 
must be order continuous.\par
Now let us assume that $X$ is \ac{msluacs} and $E^*$ has the Kadets-Klee* property and weak*-sequentially compact
unit ball. To show that $E(X)$ is \ac{msluacs} we will proceed in an analogous way to the proof of \cite{hardtke}*{Proposition 4.7},
which in turn uses techniques from the proof of \cite{dowling}*{Proposition 4}.\par
So let us take two sequences $(f_n)_{n\in \N}$, $(g_n)_{n\in \N}$ in $S_{E(X)}$ and $f\in S_{E(X)}$ such that
$\norm{f_n+g_n-2f}_{E(X)}\to 0$. Also, take a sequence $(l_n)_{n\in \N}$ of norm-one funcionals on $E(X)$ such that 
$l_n(f_n)\to 1$. As usual, $l_n$ will be represented by $[F_n]\in E^{\prime}(X^*,w^*)$ and we conclude
\begin{equation}\label{eq:3.97}
\lim_{n\to \infty}\int_S\norm{F_n(t)}\norm{f_n(t)}\,\mathrm{d}\mu(t)=1
\end{equation}
and, after passing to an appropriate subsequence,
\begin{equation}\label{eq:3.98}
\lim_{n\to \infty}\paren*{\norm{F_n(t)}\norm{f_n(t)}-F_n(t)(f_n(t))}=0 \ \ \text{a.\,e.}
\end{equation}
We also have
\begin{align*}
&\norm{2\norm{f(\cdot)}-\norm{f_n(\cdot)+g_n(\cdot)}}_E=\norm{\abs*{2\norm{f(\cdot)}-\norm{f_n(\cdot)+g_n(\cdot)}}}_E \\
&\leq\norm{\norm{2f(\cdot)-f_n(\cdot)-g_n(\cdot)}}_E=\norm{2f-f_n-g_n}_{E(X)},
\end{align*}
hence
\begin{equation}\label{eq:3.99}
\norm{2\norm{f(\cdot)}-\norm{f_n(\cdot)+g_n(\cdot)}}_E\to 0.
\end{equation}
As before we can also show
\begin{equation}\label{eq:3.100}
\norm{\norm{f_n(\cdot)}+\norm{g_n(\cdot)}}_E\to 2.
\end{equation}
Also, because of $\norm{f_n+g_n-2f}_{E(X)}\to 0$ we may pass to a further subsequence such that
\begin{equation}\label{eq:3.101}
\lim_{n\to \infty}\norm{f_n(t)+g_n(t)-2f(t)}=0 \ \ \text{a.\,e.}
\end{equation}
Let us define for every $n\in \N$
\begin{align*}
&a_n(t):=2\norm{f(t)}-\frac{1}{2}\paren*{\norm{f_n(t)}+\norm{g_n(t)}}, \\
&b_n(t):=\norm{f(t)}-\frac{1}{2}\norm{f_n(t)+g_n(t)}.
\end{align*}
Note that 
\begin{equation*}
\norm{f(t)}\leq b_n(t)+\frac{1}{2}(\norm{f_n(t)}+\norm{g_n(t)}).
\end{equation*}
So if $a_n(t)\geq 0$, then
\begin{equation*}
\abs{a_n(t)}=2\norm{f(t)}-\frac{1}{2}\paren*{\norm{f_n(t)}+\norm{g_n(t)}}\leq2\abs{b_n(t)}+\frac{1}{2}\paren*{\norm{f_n(t)}+\norm{g_n(t)}}.
\end{equation*}
If $a_n(t)<0$, then
\begin{equation*}
\abs{a_n(t)}=\frac{1}{2}\paren*{\norm{f_n(t)}+\norm{g_n(t)}}-2\norm{f(t)}\leq2\abs{b_n(t)}+\frac{1}{2}\paren*{\norm{f_n(t)}+\norm{g_n(t)}}.
\end{equation*}
So we always have
\begin{equation*}
\abs{a_n(t)}\leq 2\abs{b_n(t)}+\frac{1}{2}\paren*{\norm{f_n(t)}+\norm{g_n(t)}}.
\end{equation*}
It follows that
\begin{align*}
&\frac{1}{2}\norm{\norm{f_n(\cdot)}+\norm{g_n(\cdot)}}_E+2\norm{b_n}_E\geq\norm*{2\abs*{b_n}+\frac{1}{2}(\norm{f_n(\cdot)}+\norm{g_n(\cdot)})}_E \\
&\geq\norm{a_n}_E\geq 2-\frac{1}{2}\norm{\norm{f_n(\cdot)}+\norm{g_n(\cdot)}}_E
\end{align*}
and we can conclude with \eqref{eq:3.99} and \eqref{eq:3.100} that $\norm{a_n}_E\to 1$.\par
Using this together with \eqref{eq:3.99}, $\norm{f_n(\cdot)}+\norm{g_n(\cdot)}+2a_n=4\norm{f(\cdot)}$ and the fact
that $E$ is \ac{MLUR} we get that
\begin{equation}\label{eq:3.102}
\lim_{n\to \infty}\norm{2\norm{f(\cdot)}-\norm{f_n(\cdot)}-\norm{g_n(\cdot)}}_E=0.
\end{equation}
Again, since $E$ is \ac{MLUR} this implies
\begin{equation}\label{eq:3.103}
\lim_{n\to \infty}\norm{\norm{f_n(\cdot)}-\norm{g_n(\cdot)}}_E=0.
\end{equation}
Because of \eqref{eq:3.102} and \eqref{eq:3.103} we can pass to a further subsequence such that
\begin{equation}\label{eq:3.104}
\lim_{n\to \infty}\norm{f_n(t)}=\norm{f(t)}=\lim_{n\to \infty}\norm{g_n(t)} \ \ \text{a.\,e.}
\end{equation}
Since $B_{E^*}$ is weak*-sequentially compact we may also asssume that $(\norm{F_n(\cdot)})_{n\in \N}$ weak*-converges
to some $g\in B_{E^{\prime}}$.\par
\eqref{eq:3.102} and \eqref{eq:3.103} imply $\norm{\norm{f_n(\cdot)}-\norm{f(\cdot)}}_E\to 0$. Together with \eqref{eq:3.97}
this gives us
\begin{equation}\label{eq:3.105}
\lim_{n\to \infty}\int_S\norm{F_n(t)}\norm{f(t)}\,\mathrm{d}\mu(t)=1,
\end{equation}
hence we also have
\begin{equation*}
\int_S\norm{g(t)}\norm{f(t)}\,\mathrm{d}\mu(t)=1,
\end{equation*}
thus $\norm{g}_{E^{\prime}}=1$. Since $E^*$ has the Kadets-Klee* property it follows that $\norm{\norm{F_n(\cdot)}-g}_{E^{\prime}}\to 0$,
so if we pass again to a subsequence we may assume 
\begin{equation}\label{eq:3.106}
\lim_{n\to \infty}\norm{F_n(t)}=g(t) \ \ \text{a.\,e.}
\end{equation}
Now if we combine \eqref{eq:3.98}, \eqref{eq:3.99}, \eqref{eq:3.104} and \eqref{eq:3.106} we obtain
\begin{equation*}
\lim_{n\to \infty}\paren*{\norm{F_n(t)}\norm{f(t)}-F_n(t)(f(t))}=0 \ \ \text{a.\,e.},
\end{equation*}
since $X$ is \ac{msluacs}.\par
Using our usual argument via equi-integrability and Vitali's Lemma this leads to
\begin{equation*}
\lim_{n\to \infty}\int_A\paren*{\norm{F_n(t)}\norm{f(t)}-F_n(t)(f(t))}\,\mathrm{d}\mu(t)=0
\end{equation*}
for every $A\in \mathcal{A}$ with $\mu(A)<\infty$.\par
By the order continuity of $E$ we can derive from this
\begin{equation}\label{eq:3.107}
\lim_{n\to \infty}\int_S\paren*{\norm{F_n(t)}\norm{f(t)}-F_n(t)(f(t))}\,\mathrm{d}\mu(t)=0
\end{equation}
also in the $\sigma$-finite case (cf. the proof of Theorem \ref{thm:sluacs2}).\par
Combining \eqref{eq:3.107} and \eqref{eq:3.105} gives us $l_n(f)\to 1$ and we are done.\par
The statement about \ac{mluacs} spaces can be proved similarly.
\end{Proof}

We remark that the results proved in this section especially apply to $L^p$ spaces for $1<p<\infty$
(as we said before, for the properties \ac{acs}/\ac{luacs}/\ac{uacs} this was already proved by 
Sirotkin in \cite{sirotkin}).
\begin{corollary}\label{cor:Lp}
If $X$ is \ac{acs}/\ac{luacs}/$\text{luacs}^+$/\ac{sluacs}/$\text{sluacs}^+$/\ac{mluacs}/\ac{msluacs}/\ \ac{wuacs}/\ac{uacs} 
then for any complete, $\sigma$-finite measure space $(S,\mathcal{A},\mu)$ and any $1<p<\infty$ the 
Lebesgue-Bochner $L^p(\mu)(X)$ has the same property.
\end{corollary}
In the last section we will establish some further connections between the various properties that
we considered in this paper.

\section{Miscellaneous}\label{sec:misc}
In \cite{lovaglia} A. Lovaglia called a Banach space $X$ weakly locally uniformly rotund if for every 
sequence $(x_n)_{n\in \N}$ in $S_X$, every $x\in S_X$ and each $x^*\in S_{X^*}$ the implication 
\begin{equation*}
\norm{x_n+x}\to 2 \ \ \mathrm{and} \ \ x^*(x)=1 \ \ \Rightarrow \ \ x^*(x_n)\to 1
\end{equation*}
holds. Since this notion of weak local uniform rotundity is strictly weaker than the notion of \ac{WLUR} 
spaces that is nowadays commonly used, we will call such spaces \ac{WLUR} in the sense of Lovaglia.\footnote{A dual
Banach space will be called \ac{WLUR}* in the sense of Lovaglia if it fulfils Lovaglia's definition for all
evaluation functionals.}  By definition, a Banach space is $\mathrm{luacs}^+$ if and only if it is \ac{luacs} 
and \ac{WLUR} in the sense of Lovaglia. Also, the following is valid.
\begin{proposition}\label{prop:char luacs+}
A Banach space $X$ is $\text{luacs}^+$ if and only if $X$ is \ac{WLUR} in the sense of Lovaglia
and for all $x^*, y^*\in S_{X^*}$ with $\norm{x^*+y^*}=2$ and every $x\in S_X$ with $x^*(x)=1$ 
one also has $y^*(x)=1$.
\end{proposition}

\begin{Proof}
The necessity is clear because of \cite{hardtke}*{Proposition 2.16 (i)}. For the sufficiency we only have to 
prove that $X$ is luacs, so let us take a sequence $(x_n)_{n\in \N}$ in $S_X$ and $x\in S_X$ such that 
$\norm{x_n+x}\to 2$ as well as $x^*\in S_{X^*}$ with $x^*(x_n)\to 1$. Since $B_{X^{**}}$ is weak*-compact 
we can find $x^{**}\in B_{X^{**}}$ and a subnet $(x_{\varphi(i)})_{i\in I}$ which is weak*-convergent to $x^{**}$.
It follows that $x^{**}(x^*)=1=\norm{x^{**}}$.\par
Now fix a sequence $(y_n^*)_{n\in \N}$ in $S_{X^*}$ such that $y_n^*(x_n+x)\to 2$. Then $y_n^*(x_n)\to 1$ and 
$y_n^*(x)\to 1$. There is $y^*\in B_{X^*}$ and a subnet $(y_{\psi(j)}^*)_{j\in J}$ which is weak*-convergent 
to $y^*$. It follows that $y^*(x)=1=\norm{y^*}$. Since $X$ is \ac{WLUR} in the sense of Lovaglia we conclude
$y^*(x_n)\to 1$. It follows that $x^{**}(y^*)=1=x^{**}(x^*)$, hence $\norm{x^*+y^*}=2$.\par
Becuase of $y^*(x)=1$ our assumption imlies $x^*(x)=1$ and we are done.
\end{Proof}

The following assertion is also easy to prove (we omit the details).
\begin{proposition}\label{prop:sufficient sluacs}
If $X$ is a Banach space which \ac{WLUR} in the sense of Lovaglia and such that
$X^*$ is \ac{WLUR}* in the sense of Lovaglia then $X$ is \ac{sluacs}.
\end{proposition}

Under additional assumptions on the space $X$ it is possible to prove some more results.
\begin{proposition}\label{prop:reflexive}
Let $X$ be a reflexive Banach space.
\begin{enumerate}[\upshape(i)]
\item If $X$ is \ac{WLUR} in the sense of Lovaglia then $X$ is $\text{luacs}^+$.
\item If $X$ is \ac{sluacs} and $\text{luacs}^+$ then $X$ is \ac{wuacs}.
\item If $X$ is \ac{wuacs} and \ac{R} then $X$ is \ac{WLUR}.
\end{enumerate}
\end{proposition}

\begin{Proof}
(i) follows directly from the Proposition \ref{prop:char luacs+} and \cite{hardtke}*{Proposition 2.15}. Of the 
remaining assertions we will only prove (iii) explicitly.\par
Let $(x_n)_{n\in \N}$ be a sequence in $S_X$ and $x\in S_X$ such that $\norm{x_n+x}\to 2$. We can find a 
sequence $(x_n^*)_{n\in \N}$ in $S_{X^*}$ such that $x_n^*(x_n+x)\to 2$ and hence $x_n^*(x_n)\to 1$ and $x_n^*(x)\to 1$.\par
Since $X$ is reflexive we may assume that $(x_n^*)_{n\in \N}$ is weak*-convergent to some $y^*\in B_{X^*}$ and
$(x_n)_{n\in \N}$ is weakly convergent to some $y\in B_X$. It follows that $y^*(x)=1$ and hence $\norm{x_n^*+y^*}\to 2$.\par
Since $X$ is \ac{wuacs} the dual space $X^*$ is \ac{sluacs} (cf.\,\cite{hardtke}*{Proposition 2.16}) and 
thus (because of $x_n^*(x_n)\to 1$) we can conclude $y^*(x_n)\to 1$, whence $y^*(y)=1=y^*(x)$, which implies 
$\norm{x+y}=2$, which by the rotundity of $X$ implies $x=y$.
\end{Proof}

\begin{proposition}\label{prop:reflexive KK}
Let $X$ be a reflexive Banach space with the Kadets-Klee property.
\begin{enumerate}[\upshape(i)]
\item If $X$ is \ac{acs} then $X$ is \ac{luacs}.
\item If $X$ is \ac{WLUR} in the sense of Lovaglia then $X$ is \ac{wuacs} and $\text{sluacs}^+$.
\item If $X$ is \ac{WLUR} in the sense of Lovaglia and \ac{R} then $X$ is \ac{wuacs} and \ac{LUR}.
\end{enumerate}
\end{proposition}

\begin{Proof}
(i) Let $(x_n)_{n\in \N}$, $x$ and $y$ be as in the proof of (iii) of the previous Proposition and let $x^*\in S_{X^*}$
with $x^*(x_n)\to 1$. Then $x^*(y)=1$ and hence $\norm{y}=1$. Since $X$ has the Kadets-Klee property it follows that 
$\norm{x_n-y}\to 0$ and thus $\norm{x+y}=2$. Because $X$ is \ac{acs} we obtain $x^*(x)=1$, as desired.\par
(ii) We first show that $X$ is \ac{wuacs}. Take two sequences $(x_n)_{n\in \N}$ and $(y_n)_{n\in \N}$ in $S_X$ such 
that $\norm{x_n+y_n}\to 2$ and a functional $x^*\in S_{X^*}$ with $x^*(x_n)\to 1$. By the reflexivity of $X$ we may
assume that $(x_n)_{n\in \N}$ is weakly convergent to some $x\in B_X$. Then $x^*(x)=1$, hence $\norm{x}=1$.\par
But $X$ has the Kadets-Klee property, so this implies $\norm{x_n-x}\to 0$.\par
Now fix a sequence $(y_n^*)_{n\in \N}$ in $S_{X^*}$ such that $y_n^*(x_n)\to 1$ and $y_n^*(y_n)\to 1$. It follows 
that $y_n^*(x)\to 1$ and consequently $\norm{y_n+x}\to 2$.\par
Since $x^*(x)=1$ and $X$ is \ac{WLUR} in the sense of Lovaglia we get $x^*(y_n)\to 1$, proving that $X$ is \ac{wuacs}.\par
Now we will show that $X$ is \ac{sluacs}. Take $(x_n)_{n\in \N}$ and $x$ in $S_X$ with $\norm{x_n+x}\to 2$ and a sequence
$(x_n^*)_{n\in \N}$ in $S_{X^*}$ such that $x_n^*(x_n)\to 1$. Also, fix a sequence $(y_n^*)_{n\in \N}$ in $S_{X^*}$
with $y_n^*(x_n)\to 1$ and $y_n^*(x)\to 1$.\par
We may assume that $(x_n)_{n\in \N}$ is weakly convergent to some $y\in B_X$ and $(y_n^*)_{n\in \N}$ is weak*-convergent
to some $y^*\in B_{X^*}$. It follows that $y^*(x)=1$ and hence $\norm{y^*+y_n^*}\to 2$.\par
Since $X$ is \ac{wuacs} $X^*$ is \ac{sluacs} and thus we get $y^*(x_n)\to 1$. It follows that $y^*(y)=1$, hence $\norm{y}=1$
and $\norm{x+y}=2$. The Kadets-Klee property of $X$ gives us $\norm{x_n-y}\to 0$.\par
Because of $x_n^*(x_n)\to 1$ we can now infer $x_n^*(y)\to 1$. Since $X$ is in particular \ac{acs} this implies
$x_n^*(x)\to 1$ (cf.\,\cite{hardtke}*{Proposition 2.19}).\par
We will skip the last part of the proof, the reverse implication in the definition of $\mathrm{sluacs}^+$.\par
(iii) By (ii) $X$ is \ac{wuacs} and $\mathrm{sluacs}^+$. Let us take a sequence $(x_n)_{n\in \N}$ in $S_X$
and an element $x\in S_X$ such that $\norm{x_n+x}\to 2$. Fix a sequence $(x_n^*)_{n\in \N}$ in $S_{X^*}$ such 
that $x_n^*(x_n)=1$ for every $n\in \N$. Since $X$ is \ac{sluacs} it follows that $x_n^*(x)\to 1$.\par
Assume that $(x_n)_{n\in \N}$ is weakly convergent to $y\in B_X$ and that $(x_n^*)_{n\in \N}$ is weak*-convergent to 
$x^*\in B_{X^*}$. It follows that $x^*(x)=1$ and hence $x^*\in S_{X^*}$. Moreover, since $X$ is \ac{WLUR} in the sense of
Lovaglia we get that $x^*(x_n)\to 1$.\par
Since $(x_n)_{n\in \N}$ converges weakly to $y$ this implies $x^*(y)=1$ and hence $\norm{y}=1$. Now the
Kadets-Klee property of $X$ allows us to conclude $\norm{x_n-y}\to 0$.\par
Because of $x^*(x)=x^*(y)=1$ we must have $\norm{x+y}=2$ and thus the rotundity of $X$ implies $x=y$.
\end{Proof}

\begin{proposition}\label{prop:KK*}
Let $X$ be a Banach space such that $X^*$ has the Kadets-Klee* property and $B_{X^*}$ is weak*-sequentially compact.
\begin{enumerate}[\upshape(i)]
\item If $X$ is \ac{S} then it is also \ac{WLUR} in the sense of Lovaglia.
\item If $X^*$ is \ac{acs} then $X$ is $\text{luacs}^+$ and for all sequences $(x_n)_{n\in \N}$ in $S_X$,
$(x_n^*)_{n\in \N}$ in $S_{X^*}$ and every $x\in S_X$ with $\norm{x_n+x}\to 2$ and $x_n^*(x)\to 1$ one has
$x_n^*(x_n)\to 1$.
\item If $X^*$ is \ac{WLUR}* in the sense of Lovaglia then $X$ is \ac{sluacs}.
\end{enumerate}
\end{proposition}

\begin{Proof}
We will only prove (iii), so let $(x_n)_{n\in \N}$ and $x$ be in $S_X$ with $\norm{x_n+x}\to 2$ and $(x_n^*)_{n\in \N}$
a sequence in $S_{X^*}$ such that $x_n^*(x_n)\to 1$. Let $(y_n^*)_{n\in \N}$ be a sequence in $S_{X^*}$ with 
$y_n^*(x_n)\to 1$ and $y_n^*(x)\to 1$.\par
By assumption, we may suppose that $(y_n^*)_{n\in \N}$ is weak*-convergent to some $y^*\in B_{X^*}$. Then $y^*(x)=1$,
hence $y^*\in S_{X^*}$. By the Kadets-Klee* property of $X^*$ we must have $\norm{y_n^*-y^*}\to 0$.\par
It follows that $y^*(x_n)\to 1$, hence $\norm{x_n^*+y^*}\to 2$. Since $X^*$ is \ac{WLUR}* in the sense of Lovaglia
we obtain $x_n^*(x)\to 1$.
\end{Proof}

\address
\email

\end{document}